\RequirePackage{fix-cm}
\documentclass[smallextended]{svjour3}       
\smartqed  
\usepackage{graphicx}

\usepackage{color}
\newcommand{\CorrR}[1]{{\color{black} #1}}
\newcommand{\CorrB}[1]{{\color{black}#1}}
\newcommand{\CorrRR}[1]{{\color{black} #1}}

\usepackage{amsfonts}

\newcommand{\cvd}{\hfill{$\sqcap\!\!\!\!\sqcup$}}

\newtheorem{Theorem}{Theorem}

\newtheorem{Proposition}[Theorem]{Proposition}

\newtheorem{Definition}[Theorem]{Definition}
\newtheorem{Remark}[Theorem]{Remark}

\begin{document}

\title{A differential game with exit costs
}


\author{Fabio Bagagiolo \and Rosario Maggistro \and Marta Zoppello}


\institute{F. Bagagiolo \at
              Department of Mathematics, University of Trento, Italy  \email{fabio.bagagiolo@unitn.it}    \CorrB{\and
           R. Maggistro \at
              Department of Management, Ca' Foscari University of Venice, Italy \email{rosario.maggistro@unive.it} \and M. Zoppello \at Department of Mathematical Sciences, Politecnico di Torino, Italy \email{marta.zoppello@polito.it}
}}

\date{Received: date / Accepted: date}

\maketitle

\begin{abstract}
We study a differential game where two players separately control their own dynamics, pay a running cost, and moreover pay an exit cost (quitting the game) when they leave a fixed domain. In particular, each player has its own domain and the exit cost consists of three different exit costs, depending whether either the first player only leaves its domain, or the second player only leaves its domain, or they both simultaneously leave their own domain. We prove that, under suitable hypotheses, the lower and upper value are continuous and are, respectively, the unique viscosity solution of a suitable Dirichlet problem for a Hamilton-Jacobi-Isaacs equation. The continuity of the values relies on the existence of suitable non-anticipating strategies respecting the domain-constraint. This problem is also treated in this work. 

\end{abstract}

\section{Introduction}
\label{intro}
In some quite recent authors' works (see for example Bagagiolo \cite{bag}, Bagagiolo-Danieli \cite{bagdan}, Bagagiolo-Maggistro \cite{bagmag} and the references therein) some optimal control problems for systems governed by thermostatic dynamics are studied in the framework of dynamic programming methods and viscosity solutions of Hamilton-Jacobi-Bellman (HJB) equations. By system governed by thermostatic dynamics, here we mean  an equation as $y'=f(y,w,\alpha)$ where $\alpha$ is the measurable control and $w\in\{-1,1\}$ is the output of a hysteretical thermostat (delayed relay) subject to the evolution of some suitable components of the state-variable $y$. Hence, the switching evolution of the parameter $w$ is not directly at disposal of the external controller, but it follows some internal switching rules which are intrinsic to the system. In \cite{bag,bagdan}, the value function is proven to be the unique viscosity solution of a suitably coupled system of HJB equations, where the coupling is given by the boundary conditions in the regions where the thermostat certainly assumes a constant value (cannot switch). This is done by splitting the optimal control problem in some problems of exit time kind: in every space-region where the thermostat is constant, the problem is equivalent to an exit-time problem with unknown exit-cost given by the value function itself evaluated in the other region of constancy for $w$. Then, an ad hoc fixed point procedure is applied. Hence, a crucial starting point for such a procedure is a good theory for exit-time/exit-costs optimal control problems, in particular for what concerns the identification of the value function as the unique viscosity solution of a Dirichlet problem for HJB equations with suitable boundary conditions. Fortunately, such a good theory was quite already at disposal. 

In \cite{bag,bagdan,bagmag} some motivations and applications for studying optimal control problems with thermostatic dynamics are given. Similar motivations certainly suggest the study of differential games with thermostatic dynamics. Just think to a pursuit evasion game between two cars with automatic gears, where the switching variable(s) $w$ may represent the position of the gears. We also point out that such thermostatic dynamics is a special case of a so-called hybrid dynamics, and we refer to Gromov-Gromova \cite{grogro}, for a recent study of hybrid differential games in the framework of necessary optimality conditions. 

In order to apply to differential games some similar procedures as the ones applied to HJB for thermostatic optimal control problems, we need a good theory for exit time differential games.
Unfortunately, for differential games, the situation is rather different than from optimal control: differential games problem with exit-time and exit-costs are not so well studied in the framework of viscosity solutions for Hamilton-Jacobi equations.  Hence, before studying differential games for systems with thermostatic dynamics we need first, at least, a uniqueness results for some suitable Dirichlet problems for Hamilton-Jacobi-Isaacs (HJI) equations, in order to be able to identify the (upper and lower) values of the exit-time differential game as the unique viscosity solutions. Up to the knowledge of the authors, this paper may represent the first attempt of studying an exit-time differential games in the framework of  viscosity solutions theory for Isaacs equations with boundary conditions in the viscosity sense. The possible applications to a differential games with thermostatic dynamics will be the subject of a future work. However, we would like to point out that the differential games with exit-time and exit-costs are interesting by themselves, and not only for applications to thermostatic dynamics. To this purpose, we are going to give an example at the end of the next section.

The studied problem is the following one. We are given the controlled system

\begin{equation}
\label{eq:systems}
\begin{array}{l}
\displaystyle
\left\{
\begin{array}{ll}
\displaystyle
X'(t)=f\Big(X(t),\alpha(t)\Big),&t>0,\\
\displaystyle
X(0)=x\in\overline\Omega_X,
\end{array}
\right.\\
\\
\displaystyle
\left\{
\begin{array}{ll}
\displaystyle
Y'(t)=g\Big(Y(t),\beta(t)\Big),&t>0,\\
\displaystyle
Y(0)=y\in\overline\Omega_Y,
\end{array}
\right.
\end{array}
\end{equation}

\noindent
with $\alpha$ and $\beta$ controls, $X$ and $Y$ state-variables for the two players, respectively \CorrR{(i.e., for example, the time dependent function $t\mapsto X(t)$ is the solution (trajectory) of the first system, when the initial point $x$ and the measurable control $\alpha$ are fixed)}. The closed  sets $\overline\Omega_X$ and $\overline\Omega_Y$ (closure of the open ones $\Omega_X$ and $\Omega_Y$) are respectively the state-space for $X$ and $Y$. Denoting by $\tau_X$ and $\tau_Y$ the first exit-time from $\overline\Omega_X$ and $\overline\Omega_Y$ respectively, we define $\tau=\min\{\tau_X,\tau_Y\}$, and consider the following integral discounted cost

\[
J(x,y,\alpha,\beta)=\int_0^{\tau}e^{-\lambda t}\ell(X(t),Y(t),\alpha(t),\beta(t))dt
+e^{-\lambda\tau}\psi(X(\tau),Y(\tau)),
\]

\noindent
where, in particular, the exit-cost $\psi:\partial(\Omega_X\times\Omega_Y)\to[0,+\infty[$ is a given function, which is not required to be continuous but only separately continuous on $\partial\Omega_X\times\Omega_Y$, $\Omega_X\times\partial\Omega_Y$ and $\partial\Omega_X\times\partial\Omega_Y$. Player $X$ wants to minimize the cost whereas player $Y$ wants to maximize. The problem has then a pursuit-evasion structure. 

The ``weak" continuity hypothesis on the exit-cost $\psi$ is assumed in order to take account of the possible application to the thermostatic case. Indeed, in that case, $\overline\Omega_X$ and $\overline\Omega_Y$ represent two regions where the thermostats (one per every player) assume constant values. Hence the first exit-time represent the first switching time for the thermostat, and so exiting from $\partial\Omega_X\times\Omega_Y$ means that only the thermostat of $X$ switches, exiting from $\Omega_X\times\partial\Omega_Y$ means that only the thermostat of $Y$ switches and finally exiting from $\partial\Omega_X\times\partial\Omega_Y$ means that both thermostats switch. In every one of such three cases the new scenario of the game (which continues to run after switching) may be completely different and not related to each other. 

There are three main points which are going to be treated in this paper: 

a) continuity of the lower and upper value, 

b) derivation of suitable boundary conditions for the corresponding HJI equations, 

c) uniqueness results for those Dirichlet problems, in the sense of viscosity solutions. 

We refer the reader to Bardi-Capuzzo Dolcetta \cite{barcap} for a comprehensive account to viscosity solutions theory and applications to optimal control problems and differential games (for differential games see also Buckdahn-Cardaliaguet-Quincampoix \cite{buccarqui})

\CorrR{{\it Point a): Sections \ref{sec:controllability}, \ref{sec:lower_upper} and \ref{sec:constrained_strategies}}.} 

\noindent
The lower and upper value are respectively defined as

\[
\begin{array}{l}
\displaystyle
\underline V(x,y)=\inf_\gamma\sup_\beta J(x,y,\gamma[\beta],\beta),\\
\displaystyle
\overline V(x,y)=\sup_\xi\inf_\alpha J(x,y,\alpha,\xi[\alpha]),
\end{array}
\]

\noindent
where $\gamma$ and $\xi$ are non-anticipating strategies in the sense of Elliot-Kalton \cite{ellkal}. They are functions from the set of controls for one player to the set of controls of the other one, which do not depend on the future behavior of the control, \CorrR{in the sense that, if two controls coincide in the time interval $[0,t]$, then their images also coincide in $[0,t]$ (see Definition \ref{def:nonant})}. 

\CorrR{{\it At this level, a subpoint a1): state-constraint non-anticipating strategies, must be treated.}}

\noindent
After a necessary suitable compatibility assumptions on the exit cost (see (\ref{eq:corner}) and Remark \ref{rmrk:corner}), the main problem concerning the continuity of the values is the existence of suitable \CorrR{``state-constraint non-anticipating strategies"}. Simplifying (see Assumption 2 for more precise details), this means that, for example,  for every non-anticipating strategy $\gamma$ for player $X$, 
there exists a non-anticipating strategy $\tilde\gamma$ for $X$ such that, for every control $\beta$ for $Y$, the control $\tilde\gamma[\beta]$ makes the evolution of $X$ remain inside $\overline\Omega_X$ as long as $\beta$ makes the evolution of $Y$ remaining inside $\overline\Omega_Y$, and the cost paid by the couple of controls $(\tilde\gamma[\beta],\beta)$ is not so different from the cost paid by the couple $(\gamma[\beta],\beta)$. The difficult here is to construct such a \CorrR{state-constraint non-anticipating} strategy $\tilde\gamma$. \CorrR{The problem of the construction of a new state-constraint control that makes the state-constraint be respected and that pays a cost not so different from the cost payed by the old control, also occurs in the case of  constrained optimal control problems}, where the trajectory is constrained to remain inside a fixed set for all the times. \CorrR{In Soner \cite{son}, starting from any control, a possible construction of a state-constraint control with such properties is given}. The main assumption is the controllability on the boundary of the set (which we also assume). \CorrR{However, in that paper, the construction of the state-constraint control is done in an ``anticipating" way, that is at time $t$ the new control is constructed taking also account of the behavior of the old control for suitable times after $t$. In particular, starting from two controls that coincide in the time interval $[0,t]$, it may happen that the two corresponding state-constraint controls actually differ in $[0,t]$ (see Remark \ref{rmrk:soner}). However, in \cite{son}, there is no need of non-anticipating properties, because the argument is an optimal control problem and not a differential game, and this is why the author was not concerned with non-anticipating behaviors}. In the case of non-anticipating strategies, a similar construction as in \cite{son} is forbidden. In the present paper, using the fact that the dynamics of the two players are decoupled with respect to the space-variables  and to the controls (see (\ref{eq:systems})), \CorrRR{and also adding a decoupled feature in the controls for the cost (see (\ref{sumrunningcosts}))}, we are able to suitably adapt Soner's construction in order to get the desired \CorrR{state-constraint} non-anticipating strategy. Such assumptions on decoupled dynamics, at the present moment, seems almost necessary in order to get this kind of results. In Section \ref{sec:constrained_strategies} we actually assume a sort of more general weak decoupling of the dynamics with respect to the controls (see (\ref{nuovaf})) \CorrR{for at most one of the two players}, and we still get the desired \CorrR{state-constraint non-anticipating strategy for the player with weakly coupled dynamics}. However, that weak decoupling seems to be not immediately suitable for the results of  Sections \ref{sec:isaacs} and \ref{sec:uniqueness}. The same problem of constructing that kind of non-anticipating strategy is also studied in Bettiol-Cardaliaguet-Quincapoix \cite{betcarqui}, where the decoupled dynamics assumption is also used, and other hypotheses on the running cost are made. Other studies on constrained trajectories and non-anticipating strategies as well as on possible relations with optimal control problems and differential games can be found in Koike \cite{koi}, Bardi-Koike-Soravia \cite{barkoisor}, Cardaliaguet-Quincampoix-Saint Pierre \cite{carquisain}, Bettiol-Bressan-Vinter \cite{betbreal,betbrevin}, Bressan-Facchi \cite{breal}, Bettiol-Facchi \cite{betfac}, Bettiol-Frankowska-Vinter \cite{betfravin} and Frankowska-Marchini-Mazzola \cite{framarmaz}.

\CorrR{{\it Point b): Section \ref{sec:isaacs}.}}

\noindent
In that section, using the Dynamic Programming Principle, we prove that $\overline V$ and $\underline V$ are viscosity solutions of the corresponding HJI equation with suitable boundary conditions in the viscosity sense. As expected, such boundary conditions are determined by the exit costs on the boundary. However, in our formulation of the differential game, we are considering different exist costs, depending on which of the two players is exiting (in a state-constraint framework: which of the two players is violating the constraint). This is an important feature of a state-constraint differential game, and hence of an exit-time differential game. Which player is in charge in order to respect the constraint? Which player must be penalized when the constraint is violated? When the dynamics are not decoupled such questions have no evident answers, they may depend on the particular model under analysis. However, even if the game has a zero-sum structure (min-max), the definition of the right players' responsibility with respect to the constraint is almost always not of that kind. In our case, the dynamics are decoupled and we have different exit costs and these facts allow to rightly assign the responsibility of exit from the constraint. The compatibility condition (\ref{eq:corner}) helps to write a coherent and useful boundary condition for HJI. It says that, on the common boundary, the exit costs for the maximizing player is not larger than the cost of the minimizing one.

\CorrR{{\it Point c): Section \ref{sec:uniqueness}.}}

\noindent
We show uniqueness of $\underline V$ and $\overline V$ as viscosity solution of the corresponding Dirichlet problem for the upper and lower HJI, respectively. This is done by a rather standard double-variable technique for proving a comparison result between sub- and super-solutions,  where the boundary conditions must be treated in a non standard way.

\section{The problem}
\label{sec:1}
Let $\Omega_X\subseteq\mathbb{R}^{n}$ and $\Omega_Y\subseteq\mathbb{R}^{m}$ be two open regular sets with $n,m$ positive integers; let $A\subset\mathbb{R}^{n'},B\subset\mathbb{R}^{m'}$ be two compact sets; let $f:\mathbb{R}^n\times A\to\mathbb{R}^n$ and $g:\mathbb{R}^m\times B\to\mathbb{R}^m$ be two regular functions (i.e. bounded,  continuous and Lipschitz continuous with respect to the state-variable (their first entry) uniformly with respect to the control (their second entry)). We consider the system (\ref{eq:systems}) where $\alpha$ and $\beta$ are, respectively the measurable controls $\alpha:[0,+\infty[\to A$, $\beta:[0,+\infty[\to B$ (i.e. $\alpha\in{\cal A}$ and $\beta\in{\cal B}$).

The player $X$ uses the measurable control $\alpha$ and governs the  state variable $X(t)\in\mathbb{R}^n$. On the other hand, the player $Y$ uses the measurable control $\beta$ and governs the state variable $Y(t)\in\mathbb{R}^m$.

We are also given of a suitably regular running cost $\ell:\mathbb{R}^n\times\mathbb{R}^m\times A\times B\to[0,+\infty[$ (i.e. bounded, continuous and Lipschitz continuous with respect to the state-variables (its first two entries) uniformly with respect to the controls (its second two entries)) and of three suitably regular exit costs (i.e. bounded and continuous)

\[
\begin{array}{l}
\displaystyle
\psi_X:\partial\Omega_X\times\overline\Omega_Y\to[0,+\infty[,\\
\displaystyle
\psi_Y:\overline\Omega_X\times\partial\Omega_Y\to[0,+\infty[,\\
\displaystyle
\psi_{XY}:\partial\Omega_X\times\partial\Omega_Y\to[0,+\infty[
\end{array}
\]

\noindent
which respectively represent the costs for the exit of $X$ only (from $\overline\Omega_X$), for the exit of $Y$ only (from $\overline\Omega_Y$) and for the simultaneous exit of $X$ and $Y$. Finally we have a discount factor $\lambda>0$. 

\CorrR{
Here we collect all such hypotheses, better specifying some of them and some notations.

\begin{equation}
\label{eq:hypotheses}
\begin{array}{ll}
\displaystyle
\Omega_X\subseteq\mathbb{R}^n,\ \Omega_Y\subseteq\mathbb{R}^m\ \mbox{have } C^2-\mbox{boundary};\\
\displaystyle
A\subseteq\mathbb{R}^{n'},\ B\subseteq\mathbb{R}^{m'}\ \mbox{are compact};\ \lambda>0\\
\displaystyle
{\cal A}=\left\{\alpha:[0,+\infty[\to A\Big|\alpha\ \mbox{is measurable}\right\};\\
\displaystyle
{\cal B}=\left\{\beta:[0,+\infty[\to B\Big|\beta\ \mbox{is measurable}\right\};\\
\displaystyle
f:\mathbb{R}^n\times A\to\mathbb{R}^n,\ (x,a)\mapsto f(x,a);\ g:\mathbb{R}^m\times B\to\mathbb{R}^m,\ (y,b)\mapsto g(y,b);\\
\displaystyle
\ell:\mathbb{R}^n\times\mathbb{R}^m\times A\times B\to[0,+\infty[,\ (x,y,a,b)\mapsto\ell(x,y,a,b);\\
\displaystyle
\displaystyle
\psi_X:\partial\Omega_X\times\overline\Omega_Y\to[0,+\infty[,\ (x,y)\mapsto\psi_X(x,y);\\
\displaystyle 
\psi_Y:\overline\Omega_X\times\partial\Omega_Y\to[0,+\infty[,\ (x,y)\mapsto\psi_Y(x,y);\\
\displaystyle
\psi_{XY}:\partial\Omega_X\times\partial\Omega_Y\to[0,+\infty[,\ (x,y)\mapsto\psi_{XY}(x,y);\\
\displaystyle
f,g,\ell,\psi_X,\psi_Y,\psi_{XY}\ \mbox{are continuous and }\ \exists\ M>0\ \mbox{such that}\ \forall (x,y,a,b)\\
\displaystyle
\|f(x,a)\|,\|g(y,b)\|,|\ell(x,y,a,b)|,|\psi_X(x,y)|,|\psi_Y(x,y)|,|\psi_{XY}(x,y)|\le M;\\
\displaystyle
\exists\ L>0\ \mbox{such that}\ \forall (x_1,a),(x_2,a),(y_1,b),(y_2,b)\\
\displaystyle
\|f(x_1,a)-f(x_2,a)\|\le L\|x_1-x_2\|,\ \|g(y_1,b)-g(y_2,b)\|\le L\|y_1-y_2\|,\\
\displaystyle
|\ell(x_1,y_1,a,b)-\ell(x_2,y_2,a,b)|\le L\|(x_1,y_1)-(x_2,y_2)\|;
\displaystyle
\end{array}
\end{equation}

\noindent
In (\ref{eq:hypotheses}), $\|\cdot\|$ stays, time by time, for the corresponding Euclidean norm. Moreover, some of the hypotheses may be relaxed, as it is quite common:  the Lipschitz continuity of the running cost $\ell$ with respect to the space variable may be relaxed to a simple uniform continuity, and the regularity of the boundaries of $\Omega_X$, $\Omega_Y$ may be relaxed to a suitable piece-wise $C^2$ regularity (see for example Bagagiolo-Bardi \cite{bagbar} for such relaxation in the context of a constrained optimal control problem). Finally, in Section \ref{sec:constrained_strategies} we are going to relax a little bit the ``decoupled" feature of the dynamics $f$ and $g$ with respect to the controls.

}

When a control $\alpha$ is fixed, we define the corresponding trajectory of the first system in (\ref{eq:systems}) as $X(\cdot;x,\alpha)$; similarly we use the notation $Y(\cdot;y,\beta)$. We define the first exit time of $X$ from $\overline \Omega_X$ as

\[
\tau_X(x,\alpha)=\inf\left\{t\ge0\Big|X(t;x,\alpha)\not\in\overline\Omega_X\right\},
\]

\noindent
and, similarly, the first exit time of $Y$ from $\overline\Omega_Y$ as

\[
\tau_Y(y,\beta)=\inf\left\{t\ge0\Big|Y(t;y,\beta)\not\in\overline\Omega_Y\right\},
\]

\noindent
with the convention $\inf\emptyset=+\infty$. 

\CorrR{In the following formulas, we use the notation $\tau=\min\{(\tau_X(x,\alpha),\tau_Y(y,\beta)\}$}.
We consider the cost functional $J$, defined on  $\overline\Omega_X\times\overline\Omega_Y\times\mathcal{A}\times\mathcal{B}$, 

\[
\begin{array}{l}
\displaystyle
J(x,y,\alpha,\beta)=\\
\displaystyle
\ \ \ \ \int_0^{\tau}e^{-\lambda t}\ell(X(t;x,\alpha),Y(t;y,\beta),\alpha(t),\beta(t))dt
+e^{-\lambda\tau}\psi(X(\tau;x,\alpha),Y(\tau;y,\beta)),
\end{array}
\]

\noindent
where

\[
\begin{array}{l}
\displaystyle
e^{-\lambda\tau}\psi(X(\tau;x,\alpha),Y(\tau;y,\beta))=\\
\displaystyle
\ \ \ \ \left\{
\begin{array}{ll}
\displaystyle
e^{-\lambda\tau}\psi_X(X(\tau;x,\alpha),Y(\tau;y,\beta))&\mbox{if } \tau=\tau_X(x,\alpha)<\tau_Y(y,\beta)\le+\infty,\\
\displaystyle
e^{-\lambda\tau}\psi_Y(X(\tau;x,\alpha),Y(\tau;y,\beta))&\mbox{if } \tau=\tau_Y(y,\beta)<\tau_X(x,\alpha)\le+\infty,\\
\displaystyle
e^{-\lambda\tau}\psi_{XY}(X(\tau;x,\alpha),Y(\tau;y,\beta))&\mbox{if } \tau=\tau_X(x,\alpha)=\tau_Y(y,\beta)<+\infty\\
\displaystyle
0&\mbox{if } \tau=\inf\{+\infty,+\infty\}=+\infty.
\end{array}
\right.
\end{array}
\]

The game consists in the fact that player $X$ wants to minimize the cost $J$ and the player $Y$ wants to maximize $J$.

\CorrR{

\begin{Definition}
\label{def:nonant}

i) Let $k,\tilde k$ be two non-negative integers, and $\cal U$ be a set of measurable functions 
$u:[0,+\infty[\to\mathbb{R}^{k}$. A map that sends any $u\in{\cal U}$ to a measurable function $\tilde u:[0,+\infty[\to\mathbb{R}^{\tilde k}$ is a ``non-anticipating tuning" if, for every $u_1,u_2\in{\mathcal U}$ and for every $t\ge0$, the following holds

\[
u_1=u_2\ {\rm a.e.\ in }\ [0,t]\ \Longrightarrow\ \tilde u_1=\tilde u_2\ {\rm a.e\ in }\ [0,t].
\]

ii) The ``non-anticipating strategies for player $X$" (respectively, for player $Y$) are the elements of the set 

\[
\begin{array}{l}
\displaystyle
\Gamma=\left\{\gamma:{\mathcal B}\to{\mathcal A},\ \beta\mapsto\gamma[\beta]\ \Big|\ \forall t\ge0,\right.\\
\displaystyle
\ \ \ \ \ \ \ \ \ \beta_1=\beta_2\ \mbox{a. e. in } [0,t]\Longrightarrow \gamma[\beta_1]=\gamma[\beta_2]\ \mbox{a. e. in } [0,t]\Big\};
\end{array}
\]

(respectively, 

\[
\begin{array}{l}
\displaystyle
\chi=\left\{\xi:{\mathcal A}\to{\mathcal B},\ \alpha\mapsto\xi[\alpha]\ \Big|\ \forall t\ge0,\right.\\
\displaystyle
\ \ \ \ \ \ \ \ \ \alpha_1=\alpha_2\ \mbox{a. e. in } [0,t]\Longrightarrow \xi[\alpha_1]=\xi[\alpha_2]\ \mbox{a. e. in } [0,t]\Big\})
\end{array}
\]
\end{Definition}

Note that a non-anticipating strategy for player $X$ is a non-anticipating tuning that sends measurable controls for player $Y$ to measurable controls for player $X$. The concept of non-anticipating tuning will be used in the next sections. The concept of non-anticipating strategies is the one introduced by Elliot-Kalton in \cite{ellkal} and it is used for defining the lower and the upper value function of the differential game, respectively as
}


\[
\begin{array}{l}
\displaystyle
\underline V(x,y)=\inf_{\gamma\in\Gamma}\sup_{\beta\in{\mathcal B}}J(x,y,\gamma[\beta],\beta),\\
\displaystyle
\overline V(x,y)=\sup_{\xi\in\chi}\inf_{\alpha\in{\mathcal A}}J(x,y,\alpha,\xi[\alpha]).
\end{array}
\]

\noindent We say that the game {\it has a value} if $\underline V(x,y)=\overline V(x,y)$ for all $(x,y)\in\overline\Omega_X\times\overline\Omega_Y$.

One of the possible interesting motivations/applications of differential games with exit cost can be seen in the so-called surge tank problem, as described in Vinter-Clark \cite{Vinter_tank}, and in Falugi-Kountouriotis-Vinter \cite{fakoulvin}. Surge tanks are flow control devices, whose purpose is to prevent flow rate fluctuations
for fluids passing from one process unit to another one.
In \cite{Vinter_tank,fakoulvin} the authors, using a method given by Dupuis-McEneaney \cite{dupmce}, regard the problem as a differential game, involving
dynamics with two players $X$ and $Y$, where the objective of the $X$-player is to keep the state within a specified safe region,
despite the best efforts of the $Y$-player to drive the state out of this region. The dynamic equations of an ideal surge tank are
\begin{equation}
\label{Dyn}
\left\{
\begin{array}{ll}
\CorrR{x_1}'=x_2\\
\CorrR{x_2}'=-\alpha+\beta\\
x_1(0)=x_{1_{0}}\\
x_2(0)=x_{2_{0}}
\end{array}
\right.
\end{equation}
where $x_1$ and $x_2$ can be identified with the volume and rate of change of volume of fluid in the tank respectively, $\alpha$ is the control which regulates the rate of change of outflow and $\beta$ is the disturbance. A possible upper game is given by
\begin{equation}
\begin{array}{ll}
\sup_{\alpha\in \mathcal{A}}\inf_{\beta\in \mathcal{B}}\left(\int_0^{\tau}|\beta(t)|^2\,dt+k\tau\right)\\
\end{array}
\end{equation}
where $\tau$ denotes the first exit time from a suitable open set. The $X$ player wants to maximize the cost (to maintain the state in the safe region), whereas $Y$ wants to minimize. Note that here the dynamics are not decoupled, however, if the disturbance enters the system in a "bounded manner", then this case can be casted in the situation assumed in Section \ref{sec:constrained_strategies}. In \cite{Vinter_tank,fakoulvin} the authors are interested in bang-bang controls and in the decomposition of the problem into a collection of one player optimal control problems.

\section{Controllability}
\label{sec:controllability}

In the next section, we are going to give some regularity results and properties of the value functions. Of course, suitable, but general, hypotheses are needed. 

First of all, we assume a controllability hypothesis on the boundaries.

{\it Assumption 1}. For every $x\in\partial\Omega_X$ there exist two constant controls $a_1,a_2\in A$ such that $f(x,a_1)$ is {\it strictly entering in $\Omega_X$} and $f(x,a_2)$ is {\it strictly entering in $\mathbb{R}^n\setminus\overline\Omega_X$}. 
\CorrR{That is, denoted by $\eta_X$ the outer unit normal to the $C^2$-boundary set $\Omega_X$, then $f(x,a_1)\cdot\eta_X(x)<0$ and $f(x,a_2)\cdot\eta_X(x)>0$.} 
Similarly, for every $y\in\partial\Omega_Y$ there exist two constant controls $b_1,b_2\in B$ such that $g(y,b_1)$ is {\it strictly entering} in $\Omega_Y$ and $g(y,b_2)$ is {\it strictly entering} in $\mathbb{R}^m\setminus\overline\Omega_Y$.

Such a controllability hypothesis is essential for having the continuity of the value functions. In particular, it is linked to the existence of suitable {\it constrained non-anticipating strategy}. Indeed, the continuity of the value functions for the exit-time case presents similar features as the case of state-constraint. When we evaluate, for instance, the difference $\underline V(x_1,y_1)-\underline V(x_2,y_2)$ we need, for instance, the possibility of driving the state $X(\cdot;\cdot,x_1)$ in such a way that it remains inside $\overline\Omega_X$ until the state $X(\cdot;\cdot,x_2)$ stays inside $\overline\Omega_X$. This must be done in a way such that the variation of the cost is controlled, but the main difficulty here is the fact that it must be done in a non-anticipating way. 

\CorrR{
\begin{Definition}
\label{def:modulus}
A modulus of continuity is an increasing and continuous function $\omega:[0,+\infty[\to[0,+\infty[$ such that $\omega(0)=0$. Given a function $u:\mathbb{R}^n\to\mathbb{R}^m$, a modulus of continuity for $u$ is a modulus of continuity $\omega$ such that

\[
\|u(x)-u(y)\|\le\omega(\|x-y\|)\ \forall\ x,y.
\]

\noindent
It is well known that the existence of a modulus of continuity for $u$ is equivalent to the fact that $u$ is uniformly continuous.

\end{Definition}
}

{\it Assumption 2}. For every $T>0$, for every $K\subseteq\overline\Omega_X\times\overline\Omega_Y$ compact, there exists \CorrR{a modulus of continuity $\mathcal{O}_{T,K}$, and: 

I) for every $(x_1,y_1),(x_2,y_2)\in K$, there exists a non-anticipating tuning $\beta\mapsto\tilde\beta$ from $\cal B$ to itself (i.e. satisfying next point i), that is Definition \ref{def:nonant}), and there exists a way to associate $\tilde\gamma\in\Gamma$ to any $\gamma\in\Gamma$, such that, for every $\beta,\beta_1,\beta_2\in{\cal B}$, $\gamma\in\Gamma$, $t\ge0$, we have}

\[
\begin{array}{l}
\displaystyle
i)\ \beta_1=\beta_2\ \mbox{a.e. in } [0,t]\ \Longrightarrow\ \tilde\beta_1=\tilde\beta_2\ \mbox{a.e. in } [0,t],\\
\displaystyle
ii)\ \tau_X(x_1,\tilde\gamma[\beta])\ge \tau_X(x_2,\gamma[\beta]),\\
\displaystyle
iii)\ \tau_Y(y_2,\tilde\beta)\ge \tau_Y(y_1,\beta),\\
\displaystyle
iv)\ \|X(\tilde\tau;x_1,\tilde\gamma[\beta])-X(\tilde\tau;x_2,\gamma[\beta])\|\CorrR{\le}\mathcal{O}_{T,K}(\|x_1-x_2\|),\\
\displaystyle
v)\ \|Y(\tilde\tau;y_1,\beta)-Y(\tilde\tau;y_2,\tilde\beta)\|\CorrR{\le}\mathcal{O}_{T,K}(\|y_1-y_2\|),\\
\displaystyle
vi)\ \left|J_{\tilde\tau}(x_1,y_1,\tilde\gamma[\tilde\beta],\beta)-J_{\tilde\tau}(x_2,y_2,\gamma[\tilde\beta],\tilde\beta)\right|\\
\displaystyle
\ \ \ \ \ \ \ \ \CorrR{\le}\mathcal{O}_{T,K}(\|(x_1,y_1)-(x_2,y_2)\|),
\end{array}
\]

\noindent
where $\tilde\tau=\min\CorrR{\{}\tau_X(x_2,\gamma[\tilde\beta]),\tau_Y(y_1,\beta),T\CorrR{\}}$, and $J_{\tilde\tau}$ is the integral of the discounted running cost up to the time $\tilde\tau$.

II) A similar condition holds reversing the roles of $X$ and $Y$, $\gamma\in\Gamma$ and $\xi\in\chi$, $\alpha\in\cal A$ and $\beta\in\cal B$.

Assumption 2 is required in order to guarantee the existence of a suitable non-anticipating strategies and then prove the continuity of the values. We only need trajectories estimates on compact sets of time because the cost is discounted (the presence of the term $e^{-\lambda t}$), see the proof of Proposition \ref{prop:continuity}. Under our hypotheses, in particular the decoupling of the dynamics, and the controllability on the boundaries,  Assumption 2 holds, as it is proven in Bettiol-Cardaliaguet-Quincampoix \cite{betcarqui}. More precisely, $i), iv)$ and $v)$ are treated in Proposition 3.1 and of $vi)$ in Proposition 2.3 of \cite{betcarqui}. However, note that in \cite{betcarqui} such estimates are given for all times (and not only for compact sets) and indeed they are in a exponential fashion, which of course implies our uniform estimates on compact sets. On the other hand, conditions $ii)$ and $iii)$ just says that the constructed trajectories do not exit before the given ones. In Section \ref{sec:constrained_strategies} we are going to give a different proof of the validity of Assumption 2, modifying, in a non-anticipating manner, the proof of Soner \cite{son} for the construction of constrained controls.

\section{The lower and the upper value functions}
\label{sec:lower_upper}

By standard calculations (see for example Bardi-Capuzzo Dolcetta \cite{barcap}), $\underline V$ and $\overline V$ satisfies the usual Dynamic Programming Principle (DPP). For example, for every $t\ge0$

\begin{equation}
\label{eq:DPP}
\begin{array}{ll}
\displaystyle
\underline V(x,y)=\\
\displaystyle
\inf_{\gamma\in\Gamma}\sup_{\beta\in{\mathcal B}}\left(\int_0^\tau e^{-\lambda s}\ell(X(s;x,\gamma[\beta]),Y(s;y,\beta),\gamma[\beta](s),\beta(s))ds\right.\\
\displaystyle
\left.+e^{-\lambda\tau}\underline V(X(\tau;x,\gamma[\beta]),Y(\tau;y,\beta))\right)
\end{array}
\end{equation}

\noindent
where $\tau=\min\{t,\tau_X(x,\gamma[\beta]),\tau_Y(y,\beta)\}$.

We now assume that (see Remark \ref{rmrk:corner} for comments on it)

\begin{equation}
\label{eq:corner}
\psi_Y(x,y)\le\psi_{XY}(x,y)\le\psi_X(x,y)\ \ \forall(x,y)\in\partial\Omega_X\times\partial\Omega_Y.
\end{equation}

\begin{Proposition}
\label{prop:continuity}
Given Assumption 1, Assumption 2, hypothesis (\ref{eq:corner}) and hypotheses (\ref{eq:hypotheses}), the value functions are continuous in $\overline\Omega_X\times\overline\Omega_Y$.
\end{Proposition}

{\it Proof}. We only prove the continuity of the lower value $\underline V$. \CorrR{In particular, we are going to prove its uniform continuity in every compact set. We proceed by some steps.

1) Lets us fix $\varepsilon>0$ and take $T>0$ such that, for every trajectories and controls entering the costs, $\int_T^{+\infty} e^{-\lambda t}\ell dt+e^{-\lambda T}\psi\le\varepsilon$, where $\psi$ is any one of the exit costs $\psi_X,\psi_Y,\psi_{XY}$. This is possible, independently from the trajectories and controls inside the costs, because of the boundedness  hypotheses (\ref{eq:hypotheses}), and by the fact that the cost is discounted, i. e. $\lambda>0$. 

2) Let $K\subseteq\overline\Omega_X\times\overline\Omega_Y$ be a compact set (where we are going to prove the uniform continuity). Take another compact set $K'$, with $K\subseteq K'\subseteq\overline\Omega_X\times\overline\Omega_Y$, such that all the trajectories starting from points of $K$ belong to $K'$, for times not greater than $T$ before they possibly exit from $\overline\Omega_X$ and $\overline\Omega_Y$ respectively. That is, for example, $X(t;x,\alpha)\in K'$ for all $t\in[0,\min\{\tau_X(x,\alpha),T\}]$. Such a compact set exists by the hypotheses (\ref{eq:hypotheses}).}

\CorrR{3) 
By the compactness of $K'$ and by the continuity regularities (\ref{eq:hypotheses}), there exist $\zeta>0$ and a modulus of continuity $\omega$ such that, for every $(x,y)\in\left(\partial\Omega_X\times\partial\Omega_Y\right)\cap K'$, taking the constant controls $a_2\in A$ and $b_2\in B$ as in the Assumption 1 with respect to $x\in\partial\Omega_X$ and $y\in\partial\Omega_Y$ respectively, we have the following:  

- for every $x'\in\overline\Omega_X$ with $\|x-x'\|\le\zeta$, the trajectory starting from $x'$ with constant control $a_2$ exits from $\overline\Omega_X$ in a time interval whose length is less than $\omega(\|x-x'\|)$; for every $y'\in\overline\Omega_Y$ with $\|y-y'\|\le\zeta$, the trajectory starting from $y'$ with constant control $b_2$ exits from $\overline\Omega_Y$ in a time interval whose length is less than $\omega(\|y-y'\|)$.-

This in particular comes, besides the controllability Assumption 1 and the Lipschitz regularity of the dynamics, from the $C^2$ regularity of the boundaries, which implies that, for any piece of boundary in a compact set, the signed distance function from the boundary, $d$, is $C^2$ in a neighborhood of it. Hence one can argue estimating the signed of the composed function $t\mapsto d(z(t))$ where $z$ is  the considered trajectory (see Bardi-Capuzzo Dolcetta \cite{barcap} page 272, for a similar treatment of that function).

4) Let ${\cal O}_{T,K}$ as in Assumption 2 with respect to $T$ and $K$ fixed above, and $\overline\omega$ be a modulus of continuity for $\psi_X,\psi_Y,\psi_{XY}$ in their domains inside $K'$. We define, the following modulus of continuity

\[
\tilde\omega(r)=\max\{{\cal O}_{T,K}(r),(M+1)\omega({\cal O}_{T,K}(r)),\overline\omega((M+1)\omega({\cal O}_{T,K}(r)))\}
\]

\noindent
where $M$ is the bound of $f$, $g$, $\ell$ and of the exits costs as in (\ref{eq:hypotheses}), and $\omega$ is as in point 3).

5) Our goal is to show that, there exists $\delta>0$ such that, for all $(x_1,y_1),(x_2,y_2)\in K$ with $\|(x_1,y_1)-(x_2,y_2)\|\le\delta$, we have $|\underline V(x_1,y_1)-\underline V(x_2,y_2)|\le4\varepsilon$.

6) In the sequel, we can be concerned with the behavior of the trajectories in the time interval $[0,T]$ only. In particular, all the exit time we are going to consider will be assumed to be less than $T$. Indeed, by the previous point 1), we are going to perform the comparison analysis of the costs up to the time $T$, because, even if the game runs after $T$, then all the costs (the integrated one as well as exit ones) that will be accumulated after that time, will differ for a quantity not greater than $2\varepsilon$. Also note that, when the game run up to the time $T$, the estimate of the difference of the accumulated running costs (the integrated ones only) is standard as in the infinite-horizon case (no exit time), see for example Bardi-Capuzzo Dolcetta \cite{barcap}, Chapter VIII Proposition 1.8 and Chapter III Proposition 2.1. Hence, in the following points, and in particular in the next formulas (\ref{eq:betasegnato}) and (\ref{eq:gammasegnato}), we will assume, respectively 

\[
0\le\tau_Y(y_1,\beta)\le T,\ \ 0\le\tau_X(x_2,\gamma[\beta])\le T
\]

7) Now, take $\delta>0$ such that $\tilde\omega(\delta)\le\zeta$ and take
two arbitrary points $(x_1,y_1),(x_2,y_2)\in K$ such that $\|(x_1,y_1)-(x_2,y_2)\|\le\delta$, where $\zeta$ is given in point 3). In the sequel, $\tilde\gamma$, $\tilde\beta$ are the ones defined in Assumption 2.

8) Now, we exhibit a suitable non-anticipating tuning (see Definition \ref{def:nonant}) $\beta\mapsto\overline\beta$. For every $\beta\in\mathcal{B}$ let us define $\overline\beta\in\mathcal{B}$ for $t\in[0,T]$ as

\begin{equation}
\label{eq:betasegnato}
\overline\beta(t)=\left\{
\begin{array}{ll}
\displaystyle
\tilde\beta(t)&\mbox{if } 0\le t\le\tau_Y(y_1,\beta),\\
\displaystyle
b_2&\mbox{otherwise},
\end{array}
\right.
\end{equation}

\noindent
where $b_2\in B$ is as in Assumption 1 with respect to $Y(\tau_Y(y_1,\beta);y_1,\beta)\in\partial\Omega_Y$. Hence, by points 3)--7) and Assumption 2 (points iii) and v)), if we use the control $\overline\beta$ starting from $y_2$, then the trajectory exits from $\overline\Omega_Y$ with exit time $\tau_Y(y_2,\overline\beta)$ satisfying 

\begin{equation}
\label{eq:difftauY}
0\le\tau_Y(y_2,\overline\beta)-\tau_Y(y_1,\beta)\le\tilde\omega(\delta).
\end{equation}

In particular, $\tau_Y(y_2,\overline\beta)$ is not less than $\tau_Y(y_1,\beta)$ because, up to the time $\tau_X(y_1,\beta)$, starting from $y_2$ and using the control $\tilde\beta$, the trajectory does not exit from $\overline\Omega_Y$ (point iii) of Assumption 2). Moreover,
by point v) of Assumption 2, $Y(\tau_Y(y_1,\beta);y_2,\tilde\beta)$ is sufficiently close to the boundary (because $Y(\tau_Y(y_1,\beta);y_1,\beta)\in\partial\Omega_Y$), and hence, by point 3), the trajectory "rapidly" exits from $\overline\Omega_Y$ using the constant control $b_2$, that is (\ref{eq:difftauY}).

Also note that, being the dynamics bounded by $M$, it is

\[
\|Y(\tau_Y(y_2,\overline\beta);y_2,\overline\beta)-Y(\tau_Y(y_1,\beta);y_1,\beta)\|\le\tilde\omega(\delta).
\]

Finally, such a construction of $\overline \beta$ is a non-anticipating tuning. Indeed, if $\beta_1=\beta_2$ a. e. in the time interval $[0,T]$, then the controls $\tilde\beta_1$ and $\tilde\beta_2$ are also equal a.e. in $[0,T]$ by point i) of Assumption 2. Moreover, the trajectory $Y(\cdot;y_1,\beta_1)$ and $Y(\cdot;y_1,\beta_2)$ are also equal in the time interval $[0,t]$. Then in the interval $[0,t]$, they possibly generate the same exit time $\tau_Y(y_1,\beta_1)=\tau_Y(y_1,\beta_2)$, and hence we must have $\overline\beta_1=\overline\beta_2$ a.e. in $[0,t]$.

9) Similarly as in point 8), starting from a non-anticipating strategy for $X$ $\gamma$ we define a new non-anticipating strategy for $X$, $\overline\gamma$:

\begin{equation}
\label{eq:gammasegnato}
\overline\gamma[\beta](t)=\left\{
\begin{array}{ll}
\displaystyle
\tilde\gamma[\beta](t)&\mbox{if } 0\le t\le\tau_X(x_2,\gamma[\beta]),\\
\displaystyle
a_2&\mbox{otherwise},
\end{array}
\right.
\end{equation}

\noindent
where $a_2\in A$ is as in Assumption 1 with respect to $X(\tau_X(x_2,\gamma[\beta]);x_2,\gamma[\beta])\in\partial\Omega_X$. In this case, by Assumption 2 points ii) and iv), and by the previous points of this proof, we have, as in point 8),  

\begin{equation}
\begin{array}{ll}
\displaystyle
\label{eq:difftauX}
0\le\tau_X(x_1,\overline\gamma[\beta])-\tau_X(x_2,\gamma[\beta])\le\tilde\omega(\delta)\\
\displaystyle
\|X(\tau_X(x_1,\overline\gamma[\beta]);x_1,\overline\gamma[\beta])-X(\tau_X(x_2,\gamma[\beta]);x_2,\gamma[\beta])\|\le\tilde\omega(\delta).
\end{array}
\end{equation}

Note that, by our hypotheses, in particular because $\tilde\gamma$ is a non-anticipating strategy (the one given by Assumption 2), we have that $\overline\gamma$ is also a non-anticipating strategy. Moreover, for every $\gamma\in\Gamma$, we also consider the following non-anticipating strategy $\overline{\overline\gamma}\in\Gamma$, defined as $\overline{\overline\gamma}[\beta]=\overline\gamma[\overline\beta]$, for all $\beta\in\mathcal{B}$.

10) For suitable $\gamma_2\in\Gamma$ and $\beta_1\in{\cal B}$, by definition of infimum and supremum, we have

\[
\begin{array}{l}
\displaystyle
\underline V(x_1,y_1)-\underline V(x_2,y_2)\le\inf_{\gamma\in\Gamma}\sup_{\beta\in\mathcal{B}}J(x_1,y_1,\gamma[\beta],\beta)-\sup_{\beta\in\mathcal{B}}J(x_2,y_2,\gamma_2[\beta],\beta)+\varepsilon\\
\displaystyle
\le\sup_{\beta\in\mathcal{B}}J(x_1,y_1,\overline{\overline\gamma}_2[\beta],\beta)-\sup_{\beta\in\mathcal{B}}J(x_2,y_2,\gamma_2[\beta],\beta)+\varepsilon\\
\displaystyle
\le J(x_1,y_1,\overline{\overline\gamma}_2[\beta_1],\beta_1)-\sup_{\beta\in\mathcal{B}}J(x_2,y_2,\gamma_2[\beta],\beta)+2\varepsilon\\
\le J(x_1,y_1,\overline{\overline\gamma}_2[\beta_1],\beta_1)-J(x_2,y_2,\gamma_2[\overline\beta_1],\overline\beta_1)+2\varepsilon\\
=J(x_1,y_1,\overline\gamma_2[\overline\beta_1],\beta_1)-J(x_2,y_2,\gamma_2[\overline\beta_1],\overline\beta_1)+2\varepsilon.
\end{array}
\]

\noindent
Now, we define $\overline\tau_{12}=\min(\tau_X(x_2,\gamma_2[\overline\beta_1]),\tau_Y(y_1,\beta_1))$. Hence we have

\[
\begin{array}{l}
\displaystyle
\underline V(x_1,y_1)-\underline V(x_2,y_2)\\
\displaystyle
\le J_{\overline\tau_{12}}(x_1,y_1,\overline\gamma_2[\overline\beta_1],\beta_1)+\Psi_1-J_{\overline\tau_{12}}(x_2,y_2,\gamma_2[\overline\beta_1],\overline\beta_1)-\Psi_2+3\varepsilon,
\end{array}
\]

\noindent
where $J_{\overline\tau_{12}}$ stays for the integral cost up to the time $\overline\tau_{12}$, and $\Psi_1$ and $\Psi_2$  stay to indicate two possible sums of remaining integral cost and discounted exit costs paid by the trajectories (their explicit formulations will be given in the following points).

By our definition, note that, up to the time $\overline\tau_{12}$, it is $\overline\gamma_2[\overline\beta_1]=\tilde\gamma_2[\tilde\beta_1]$, $\gamma_2[\overline\beta_1]=\gamma_2[\tilde\beta_1]$, and $\overline\beta_1=\tilde\beta_1$. hence, by Assumption 2 point vi), and by the previous points of this proof, we have

\[
\underline V(x_1,y_1)-\underline V(x_2,y_2)\le\Psi_1-\Psi_2+\tilde\omega(\delta)+3\varepsilon.
\]

11) We now analyze $\Psi_1-\Psi_2$. We have some sub-cases.

11a) $\overline\tau_{12}=\tau_X(x_2,\gamma_2[\overline\beta_1])\le\tau_X(x_1,\overline\gamma_2[\overline\beta_1])<\tau_Y(y_1,\beta_1)$. Then

\[
\begin{array}{l}
\displaystyle
\Psi_1=\int_{\overline\tau_{12}}^{\tau_X(x_1,\overline\gamma_2[\overline\beta_1])}e^{-\lambda s}\ell ds\\
\displaystyle
\ +e^{-\lambda\tau_X(x_1,\overline\gamma_2[\overline\beta_1])}\psi_X(X(\tau_X(x_1,\overline\gamma_2[\overline\beta_1]);x_1,\overline\gamma_2[\overline\beta_1]),Y(\tau_X(x_1,\overline\gamma_2[\overline\beta_1]);y_1,\beta_1)),\\
\displaystyle
\Psi_2=e^{-\lambda\overline\tau_{12}}\psi_X(X(\overline\tau_{12};x_2,\gamma_2[\overline\beta_1]),Y(\overline\tau_{12};y_2,\overline\beta_1))
\end{array}
\]

\noindent
From this, by the previous points, the uniform continuity of $\psi_X$ in $\overline\Omega_X\cap K$, and the definition of $\tilde\omega$, we get

\[
|\Psi_1-\Psi_2|\le2\tilde\omega(\delta),
\]

11b) $\overline\tau_{12}=\tau_X(x_2,\gamma_2[\overline\beta_1])\le\tau_Y(y_1,\beta_1)\le\tau_X(x_1,\overline\gamma_2[\overline\beta_1])$. Hence $\Psi_1$ is the integral of the cost on the interval $[\overline\tau_{12},\tau_Y(y_1,\beta_1)]$ (whose length is not greater than $\tilde\omega(\delta)$, (\ref{eq:difftauX})), plus the discounted exit cost $\tilde\Psi_1$, the latter being (respectively for $\tau_Y(y_1,\beta_1)<\tau_X(x_1,\overline\gamma_2[\overline\beta_1])$ and for $\tau_Y(y_1,\beta_1)=\tau_X(x_1,\overline\gamma_2[\overline\beta_1])$):

\[
\begin{array}{ll}
\displaystyle
\tilde\Psi_1=e^{-\lambda\tau_Y(y_1,\beta_1)}\psi_Y(X(\tau_Y(y_1,\beta_1);x_1,\overline\gamma_2[\overline\beta_1]),Y(\tau_Y(y_1,\beta_1);y_1,\beta_1)),\\
\displaystyle
\tilde\Psi_1=e^{-\lambda\tau_Y(y_1,\beta_1)}\psi_{XY}(X(\tau_Y(x_1,\overline\gamma_2[\overline\beta_1]);x_1,\overline\gamma_2[\overline\beta_1]),Y(\tau_Y(y_1,\beta_1);y_1,\beta_1))
\end{array}
\]

\noindent
Similarly, $\Psi_2$, which has no integral part, is (respectively for $\overline\tau_{12}=\tau_X(x_2,\gamma_2[\overline\beta_1])<\tau_Y(y_2,\overline\beta_1)$ and for $\overline\tau_{12}=\tau_X(x_2,\gamma_2[\overline\beta_1])=\tau_Y(y_2,\overline\beta_1)$)
\[
\begin{array}{ll}
\displaystyle
\Psi_2=e^{-\lambda\overline\tau_{12}}\psi_X(X(\overline\tau_{12};x_2,\gamma_2[\overline\beta_1]),Y(\overline\tau_{12};y_2,\overline\beta_1)\\
\displaystyle
\Psi_2=e^{-\lambda\overline\tau_{12}}\psi_{XY}(X(\overline\tau_{12};x_2,\gamma_2[\overline\beta_1]),Y(\overline\tau_{12};y_2,\overline\beta_1).
\end{array}
\]

Using (\ref{eq:corner}), all the estimates of the previous points, and the uniform continuity in the compact sets of the exit costs, together with the definition of $\tilde\omega$, we have

\[
\begin{array}{ll}
\displaystyle
\tilde\Psi_1\le e^{-\lambda\overline\tau_{12}}\psi_{XY}(X(\tau_X(x_2,\gamma_2[\overline\beta_1]);x_2,\gamma_2[\overline\beta_1]),Y(\tau_Y(y_1,\beta_1);y_1,\beta_1))+\tilde\omega(\delta),\\
\displaystyle
\Psi_2\ge e^{-\lambda\overline\tau_{12}}\psi_{XY}(X(\tau_X(x_2,\gamma_2[\overline\beta_1]);x_2,\gamma_2[\overline\beta_1]),Y(\tau_Y(y_1,\beta_1);y_1,\beta_1))-\tilde\omega(\delta),
\end{array}
\]

\noindent
from which

\[
|\Psi_1-\Psi_2|\le3\tilde\omega(\delta).
\]


11c) The cases $\overline\tau_{12}=\tau_Y(y_1,\beta_1)\le\tau_Y(y_2,\overline\beta_1)<\tau_X(x_2,\gamma_2[\overline\beta_1])$ and $\overline\tau_{12}=\tau_Y(y_1,\beta_1)\le\tau_X(x_2,\gamma_2[\overline\beta_1])\le\tau_Y(y_2,\overline\beta_1)$ are similar to points 11a) and 11b) respectively.


12) Putting together the points 10) and 11), and reversing the role of $(x_1,y_1)$ and $(x_2,y_2)$, we then get point 5) taking $\delta$ sufficiently small. The proof is concluded.


\cvd

}

\begin{Remark}
\label{rmrk:corner}
Note that, whenever hypothesis (\ref{eq:corner}) is not satisfied, then the continuity is not guaranteed. For example, if $\psi_X<\psi_Y$ on $\partial\Omega_X\times\partial\Omega_Y$, then, using also the controllability hypothesis on the boundaries, we can approximate points on $\partial\Omega_X\times\partial\Omega_Y$ by points in $\partial\Omega_X\times\Omega_Y$ and by points in $\Omega_X\times\partial\Omega_Y$ where, respectively, $\underline V\le\psi_X$ and $\underline V\ge\psi_Y$. And this fact makes immediately fail the continuity of $\underline V$ on points of $\partial\Omega_X\times\partial\Omega_Y$.
\end{Remark}

\section{The lower and the upper Hamilton-Jacobi-Isaacs equations and boundary conditions}
\label{sec:isaacs}

For every $(x,y)\in\overline\Omega_X\times\overline\Omega_Y$ and for every $(p,q)\in\mathbb{R}^n\times\mathbb{R}^m$, we introduce the Upper Hamiltonian \CorrB{("$\cdot$" stays for the scalar product)}

\[
UH(x,y,p,q)=\min_{b\in B}\max_{a\in A}\left\{-f(x,a)\cdot p-g(y,b)\cdot q-\ell(x,y,a,b)\right\}
\]

\noindent
and the Lower Hamiltonian

\[
LH(x,y,p,q)=\max_{a\in A}\min_{b\in B}\left\{-f(x,a)\cdot p-g(y,b)\cdot q-\ell(x,y,a,b)\right\}.
\]

In the sequel, for a function $u$, $u_x(x,y)$ and $u_y(x,y)$ will denote the gradient with respect to $x$ and with respect to $y$, respectively.

\begin{Theorem}
\label{thm:UH}
\CorrB{Given Assumption 1, regularities (\ref{eq:hypotheses}) and condition (\ref{eq:corner})}, the lower value $\underline V$ satisfies the following (upper) problem in the viscosity sense (here expressed for a generic function $u:\overline\Omega_X\times\overline\Omega_Y\to\mathbb{R}$)

\begin{equation}
\label{eq:isaacs}
\left\{
\begin{array}{ll}
\displaystyle
\lambda u(x,y)+UH(x,y,u_x(x,y),u_y(x,y))=0&\mbox{in } \Omega_X\times\Omega_Y,\\
\displaystyle
u=\psi_X&\mbox{on } \partial\Omega_X\times\Omega_Y,\\
\displaystyle
u=\psi_Y&\mbox{on } \Omega_X\times\partial\Omega_Y,\\
\displaystyle
u=\psi_Y\ \mbox{or } u=\psi_X&\mbox{on } \partial\Omega_X\times\partial\Omega_Y.
\end{array}
\right.
\end{equation}
\end{Theorem}

By solutions in the viscosity sense we mean the following: let $\varphi\in C^1(\overline\Omega_X\times\overline\Omega_Y)$ and $(x_0,y_0)\in\overline\Omega_X\times\overline\Omega_Y$, then the following facts i) and ii) hold true:

i) if $(x_0,y_0)$ is a point of local maximum for $u-\varphi$, with respect to $\overline\Omega_X\times\overline\Omega_Y$, then we have the following four implications (one per every line)

\begin{equation}
\label{eq:subsol}
\begin{array}{l}
\left.
\begin{array}{ll}
\displaystyle
(x_0,y_0)\in\Omega_X\times\Omega_Y,\\
\displaystyle
(x_0,y_0)\in\partial\Omega_X\times\Omega_Y,\ u(x_0,y_0)>\psi_X(x_0,y_0),\\
\displaystyle
(x_0,y_0)\in\Omega_X\times\partial\Omega_Y,\ u(x_0,y_0)>\psi_Y(x_0,y_0),\\
\displaystyle
(x_0,y_0)\in\partial\Omega_X\times\partial\Omega_Y,\  \psi_X(x_0,y_0)\neq u(x_0,y_0)>\psi_Y(x_0,y_0)\\
\end{array}
\right\}\ \Longrightarrow\\
\\
\displaystyle
\ \ \ \ \lambda u(x_0,y_0)+UH(x_0,y_0,\varphi_x(x_0,y_0),\varphi_y(x_0,y_0))\le0;
\end{array}
\end{equation}

ii)  if $(x_0,y_0)$ is a point of local minimum for $u-\varphi$, with respect to $\overline\Omega_X\times\overline\Omega_Y$, then we have the following four implications (one per every line)

\begin{equation}
\label{eq:supersol}
\begin{array}{l}
\left.
\begin{array}{ll}
\displaystyle
(x_0,y_0)\in\Omega_X\times\Omega_Y,\\
\displaystyle
(x_0,y_0)\in\partial\Omega_X\times\Omega_Y,\ u(x_0,y_0)<\psi_X(x_0,y_0),\\
\displaystyle
(x_0,y_0)\in\Omega_X\times\partial\Omega_Y,\ u(x_0,y_0)<\psi_Y(x_0,y_0),\\
\displaystyle
(x_0,y_0)\in\partial\Omega_X\times\partial\Omega_Y,\  \psi_Y(x_0,y_0)\neq u(x_0,y_0)<\psi_X(x_0,y_0)\\
\end{array}
\right\}\ \Longrightarrow\\
\\
\displaystyle
\ \ \ \ \lambda u(x_0,y_0)+UH(x_0,y_0,\varphi_x(x_0,y_0),\varphi_y(x_0,y_0))\ge0;
\end{array}
\end{equation}

If $u$ satisfies i), it is said to be a subsolution; if it satisfies ii), it is said to be a supersolution. \CorrB{The equation in the first line of (\ref{eq:isaacs}) is called the (upper) Hamilton-Jacobi-Isaacs equation (or simply Isaacs). The implications given by the second, third and fourth lines of (\ref{eq:subsol})--(\ref{eq:supersol}) represent the boundary conditions in the viscosity sense.

Note that, in the formulation of (\ref{eq:isaacs}), the intermediate exit cost $\psi_{XY}$ does not play any role. We refer the reader to Remark \ref{rmrk:Nash} for more details on this fact.

\begin{Remark}
Under the same hypotheses of Theorem \ref{thm:UH}, and with the same definitions for solutions and  the boundary conditions, the upper value function $\overline V$ is a solution in the viscosity sense of

\begin{equation}
\label{eq:isaacs2}
\left\{
\begin{array}{ll}
\displaystyle
\lambda u(x,y)+LH(x,y,u_x(x,y),u_y(x,y))=0&\mbox{in } \Omega_X\times\Omega_Y,\\
\displaystyle
u=\psi_X&\mbox{on } \partial\Omega_X\times\Omega_Y,\\
\displaystyle
u=\psi_Y&\mbox{on } \Omega_X\times\partial\Omega_Y,\\
\displaystyle
u=\psi_Y\ \mbox{or } u=\psi_X&\mbox{on } \partial\Omega_X\times\partial\Omega_Y.
\end{array}
\right.
\end{equation}

\end{Remark}
}

{\it Proof of Theorem \ref{thm:UH}}. Let us note that, by our hypothesis of controllability on the boundaries, we have the following inequalities

\[
\underline V\le\psi_X\ \mbox{on } \partial\Omega_X\times\overline\Omega_Y,\ \ \underline V\ge\psi_Y\ \mbox{on } \overline\Omega_X\times\partial\Omega_Y.
\]

\noindent
\CorrB{This is because both players, from their own boundary, can immediately exit, stopping the game and paying the corresponding exit cost.}

\noindent
Hence, for the subsolution case, we only need to prove the validity of the Isaacs equation in $\Omega_X\times\Omega_Y$, the boundary condition on $\Omega_X\times\partial\Omega_Y$ and the boundary condition on $\partial\Omega_X\times\partial\Omega_Y$ for only the case $\psi_Y<\underline V<\psi_X$. In the same way, for the supersolution case, we only need to prove the Isaacs equation in $\Omega_X\times\Omega_Y$, the boundary condition on $\partial\Omega_X\times\Omega_Y$ and the boundary condition on $\partial\Omega_X\times\partial\Omega_Y$ for only the case $\psi_X>\underline V>\psi_Y$.
Since the validity of the Isaacs equation in $\Omega_X\times\Omega_Y$ is standard (see for instance Bardi-Capuzzo Dolcetta \cite{barcap} page 438), we only concentrate on the boundary conditions. 

{\it Supersolution.} Let $(x_0,y_0)\in\partial\Omega_X\times\Omega_Y$ be of minimum for $\underline V-\varphi$, and by absurd hypothesis, let us suppose that

\begin{equation}
\label{eq:absurd}
\begin{array}{l}
\displaystyle
\underline V(x_0,y_0)<\psi_X(x_0.y_0),\\
\displaystyle
\lambda\underline V(x_0,y_0)+UH(x_0,y_0,\varphi_x(x_0,y_0),\varphi_y(x_0,y_0))<0.
\end{array}
\end{equation}

\noindent
Of course, it is not restrictive to suppose that $\underline V(x_0,y_0)=\varphi(x_0,y_0)$ and that $\lambda =1$. Also using this assumption, we have that, for some $\varepsilon>0$, and for every $(x,y)\in B((x_0,y_0),\varepsilon)\cap\overline\Omega_X\times\overline\Omega_Y$ (here $B((x,y),r)$ stays for a ball of $\mathbb{R}^n\times\mathbb{R}^m$ with center in $(x,y)$ and radius $r>0$)

\begin{equation}
\label{eq:afterabsurd}
\begin{array}{l}
\varphi(x,y)+UH(x,y,\varphi_x(x,y),\varphi_y(x,y))\le-\varepsilon,\\
\displaystyle
\underline V(x,y)\ge\varphi(x,y).
\end{array}
\end{equation}

\noindent
Moreover, we can also suppose that $y\in\Omega_Y$ for all $(x,y)\in B((x_0,y_0),\varepsilon)$. Now, let $t>0$ be such that, for every $\gamma\in\Gamma$ and $\beta\in\mathcal{B}$, $(X(s;x,\gamma[\beta]),Y(s;y,\beta))\in B((x_0,y_0),\varepsilon)$ for all $(x,y)\in B((x_0,y_0),\varepsilon/2)$, for all $0\le s\le t$. Let us define $\delta=\varepsilon(1-e^{-t})/2>0$ and take $\gamma\in\Gamma$ such that, by the Dynamic Programming Principle (\ref{eq:DPP}), for any $\beta\in\mathcal{B}$ (note that it is $\tau_Y(y_0,\beta)\ge t$)

\[
\underline V(x_0,y_0)>
\left\{
\begin{array}{l}
\displaystyle
-\delta+\int_0^te^{-s}\ell(X(s;x_0,\gamma[\beta]),Y(s;y_0,\beta),\gamma[\beta](s),\beta(s))ds\\
\displaystyle
\ \ \ \ \ +e^{-t}\underline V(X(t;x_0,\gamma[\beta]),Y(t;y_0,\beta))\\
\displaystyle
\mbox{if } t<\tau_X(x_0,\gamma[\beta]),\\
\displaystyle
-\delta+\int_0^{\tau_X(x_0,\gamma[\beta])}e^{-s}\ell(X(s;x_0,\gamma[\beta]),Y(s;y_0,\beta),\gamma[\beta](s),\beta(s))ds\\
\displaystyle
\ \ \ \ \ +e^{-\tau_X(x_0,\gamma[\beta])}\psi_X(X(\tau_X(x_0,\gamma[\beta]);x_0,\gamma[\beta]),Y(\tau_X(x_0,\gamma[\beta]);y_0,\beta))\\
\displaystyle
\mbox{if } t\ge\tau_X(x_0,\gamma[\beta]).
\end{array}
\right.
\]

\noindent
In particular, we can take $\beta\equiv b$, with $b\in B$ arbitrary. 
Note that we can certainly suppose that $0<t<\tau_X(x_0,\gamma[\beta])$ for all $\beta$. Indeed, if not, for three sequences $t_n\to0^+$, $\gamma_n\in\Gamma$, and $\beta_n\in{\cal B}$ we would have $[0,t_n]\ni\tau_X(x_0,\gamma_n[\beta_n])\to0^+$ and (dropping the notations of the trajectories in the entries)

\[
-\delta_n+\int_0^{\tau_X}e^{-s}\ell ds+e^{-\tau_X}\psi_X<\underline V(x_0,y_0)\le\int_0^{\tau_X}e^{-s}\ell ds+e^{-\tau_X}\psi_X,
\]

\noindent
where $\delta_n=\varepsilon(1-e^{-t_n})/2$. But then, in the limit, this implies $\underline V(x_0,y_0)=\psi_X(x_0,y_0)$, against our absurd hypothesis (\ref{eq:absurd}). Hence we get, for every $b\in B$, (here, $X(\cdot)$ and $Y(\cdot)$ stay for the trajectories starting from $x_0$ and $y_0$ with controls $\gamma[b]$ and $b$ respectively)

\[
\begin{array}{l}
\displaystyle
0=\underline V(x_0,y_0)-\varphi(x_0,y_0)\\
\displaystyle
>-\delta+\int_0^{t}e^{-s}\ell(X(s),Y(s),\gamma[b](s),b) ds+e^{-t}\underline V(X(t),Y(t))-\varphi(x_0,y_0)\\
\displaystyle
\ge-\delta+\int_0^{t}e^{-s}\ell(X(s),Y(s),\gamma[b](s),b) ds+e^{-t}\varphi(X(t),Y(t))-\varphi(x_0,y_0)\\
\displaystyle
=-\delta+\int_0^te^{-s}\Big(-\varphi(X(s),Y(s))+\\
\displaystyle
\ \ \ \ (\varphi_x(X(s),Y(S)),\varphi_y(X(s),Y(s)))\cdot(f(X(s),\gamma[b](s)),g(Y(s),b))+\\
\displaystyle
\ \ \ \ \ \ \ \ \ell(X(s),Y(s),\gamma[b](s),b)\Big)ds\\
\displaystyle
\ge-\delta\\
\displaystyle-\int_0^te^{-s}\max_{a\in A}\left(\varphi-(\varphi_x,\varphi_y)\cdot(f(X,Y,a,b),g(X,Y,a,b))-\ell(X,Y,a,b)\right)ds.
\end{array}
\]

\noindent
In the previous formula, we used (\ref{eq:afterabsurd}), the fact that $\gamma[b](s)$ is almost everywhere an element of $A$, and the time derivative of the function $s\mapsto e^{-s}\varphi(X(s),Y(s))$. By the arbitrariness of $b\in B$, and by the definitions of the Hamiltonian $UH$ and of $\delta$, we then get the following contradiction

\[
\begin{array}{l}
\displaystyle
0>-\delta-\int_0^{t}e^{-s}\left(\varphi+UH\Big(X(s),Y(s),\varphi_x(X(s),Y(s)),\varphi_y(X(s),Y(s))\Big)\right)ds\\
\displaystyle
\ge-\delta+\int_0^te^{-s}\varepsilon ds>0.
\end{array}
\]

Let us now consider the case $(x_0,y_0)\in\partial\Omega_X\times\partial\Omega_Y$. We only have the case $\psi_Y(x_0,y_0)<\underline V(x_0,y_0)<\psi_X(x_0,y_0)$, and we can again restrict ourselves to the case $0<t<\min\{\tau_Y,\tau_X\}$, and then, arguing as before, we get the conclusion.

{\it Subsolution.} We only treat the case $(x_0,y_0)\in\Omega_X\times\partial\Omega_Y$ of maximum for $\underline V-\varphi$, the other cases being similar. Let us take $\varepsilon$, $t>0$ and $\delta>0$ in a similar way as before (changing the role of $\partial\Omega_X$ and of $\partial\Omega_Y$). For every $a\in A$ (and so for every constant strategies $\gamma\equiv a$) we find $\beta\in\mathcal{B}$ such that $0<t<\tau_Y(y_0,\beta)$ independently on $a$ (otherwise we get the contradiction $\underline V(x_0,y_0)\le\psi_Y(x_0,y_0)$) and

\[
\underline V(x_0,y_0)<\delta+\int_0^te^{-s}\ell(X(s;x_0,a),Y(s;y_0,\beta),a,\beta(s))ds+e^{-t}\underline V(X(t),Y(t)).
\]

\noindent
Hence, for every $a\in A$, we get 

\[
\begin{array}{l}
\displaystyle
0=\underline V(x_0,y_0)-\varphi(x_0,y_0)\\
\displaystyle
<\delta-\int_0^te^{-s}\min_{b\in b}\left(\varphi-(\varphi_x,\varphi_y)\cdot(f(X,Y,a,b),g(X,Y,a,b))-\ell(X,Y,a,b)\right)ds,
\end{array}
\]

\noindent 
from which the contradiction (by the arbitrariness of $a\in A$)

\[
\begin{array}{l}
\displaystyle
0<\delta-\int_0^{t}e^{-s}\left(\varphi+UH\Big(X(s),Y(s),\varphi_x(X(s),Y(s)),\varphi_y(X(s),Y(s))\Big)\right)ds\\
\displaystyle
\le+\delta-\int_0^te^{-s}\varepsilon<0.
\end{array}
\]
\cvd

\section{Uniqueness}\label{sec:uniqueness}

We are going to use the following inner-cone property of the boundaries. There exist two bounded continuous functions $\eta_X:\overline\Omega_X\to\mathbb{R}^n$, $\eta_Y:\overline\Omega_Y\to\mathbb{R}^m$, and two real positive continuous functions $c,d:\overline\Omega_X\cup\overline\Omega_Y\to]0,+\infty[$ such that, for all $x\in\overline\Omega_X$ (respectively $y\in\overline\Omega_Y$) and for all $s\in]0,d(x)]$ (respectively $s\in]0,d(y)]$), 

\begin{equation}
\label{eq:cone}
B(x+s\eta_X(x),c(x)s)\subseteq\Omega_X,\ \ (\mbox{respectively }B(y+s\eta_Y(y),c(y)s)\subseteq\Omega_Y).
\end{equation}

\CorrB{Condition (\ref{eq:cone}) roughly means that at every point there is a small cone with vertex in that point and contained, besides the vertex, in the interior of the set. Note that this is essentially a boundary regularity property, and it is certainly satisfied under the $C^2$-regularity hypothesis (\ref{eq:hypotheses}), where you can take, on the boundary, the unit interior normal as $\eta_X$ and $\eta_Y$. Finally note that in any compact subset, $\eta_X$ and $\eta_Y$ can be assumed uniformly continuous, and $c$ and $d$ just two positive constants.
}

\begin{Theorem}
\label{thm:uniqueness}
Let Assumption 1, (\ref{eq:hypotheses}) (and hence (\ref{eq:cone})), and (\ref{eq:corner}) hold. Then the lower value $\underline V$ (respectively, the upper value function $\overline V$) is the unique bounded and continuous function on $\overline\Omega_X\times\overline\Omega_Y$ which is a solution of (\ref{eq:isaacs}) (respectively, of (\ref{eq:isaacs2})) in the viscosity sense.
\end{Theorem}

\CorrB{

We are going to only prove uniqueness for (\ref{eq:isaacs}) among continuous and bounded functions, from which the theorem follows because $\underline V$ is a continuous and bounded solution in the viscosity sense. As usual, we prove such a uniqueness result by proving a comparison result between sub- and supersolutions, and we will refer to the standard double variable technique, and in particular to the ``constrained" double variable technique of Soner (see Bardi-Capuzzo Dolcetta \cite{barcap} pages 278--281) Let $u,v:\overline\Omega_X\times\overline\Omega_Y\to\mathbb{R}$ be two bounded and continuous sub- and supersolution, respectively. We are going to prove that $u\le v$ on $\overline\Omega_X\times\overline\Omega_Y$, from which the conclusion follows because every solution in the viscosity sense is simultaneously a sub-and a supersolution.

The standard procedure is as here explained. By contradiction, let us suppose

\begin{equation}
\label{eq:absurd}
\sup_{(x,y)\in\overline\Omega_X\times\overline\Omega_Y}(u(x,y)-v(x,y))=m>0,
\end{equation}

\noindent
and try to construct two test functions $\varphi_1$ and $\varphi_2$ such that, for some $(\overline x_1,\overline y_1),(\overline x_2,\overline y_2)\in\overline\Omega_X\times\overline\Omega_Y$ and some $k>0$,

\begin{equation}
\label{eq:explanation}
\begin{array}{ll}
\displaystyle
u(\overline x_1,\overline y_1)-v(\overline x_2,\overline y_2)>k\\
\displaystyle
u-\varphi_1\ \mbox{has a local maximum at } (\overline x_1,\overline y_1)\\
\displaystyle
v-\varphi_2\ \mbox{has a local minimum at } (\overline x_2,\overline y_2)\\
\displaystyle
UH(\overline x_2,\overline y_2,(\varphi_2)_x,(\varphi_2)_y)-
UH(\overline x_1,\overline y_1,(\varphi_1)_x,(\varphi_1)_y)<k\\
\displaystyle
u(\overline x_1,\overline y_1)+UH(\overline x_1,\overline y_1,(\varphi_1)_x,(\varphi_1)_y)\le0\\
\displaystyle
v(\overline x_2,\overline y_2)+UH(\overline x_2,\overline y_2,(\varphi_2)_x,(\varphi_2)_y)\ge0.
\end{array}
\end{equation}

\noindent
We then get the contradiction

\begin{equation}
\label{eq:contradiction}
\begin{array}{ll}
\displaystyle
k<u(\overline x_1,\overline y_1)-v(\overline x_2,\overline y_2)\le\\
\displaystyle
UH(\overline x_2,\overline y_2,(\varphi_2)_x,(\varphi_2)_y)-
UH(\overline x_1,\overline y_1,(\varphi_1)_x,(\varphi_1)_y)<k
\end{array}
\end{equation}

The main ingredients for this procedure are some continuity and uniform continuity properties satisfied by the Hamiltonian $UH$ (see Bardi-Capuzzo Dolcettta \cite{barcap} page 443 formula (2.1)), and the construction of a suitable penalizing function $\tilde\phi:(\overline\Omega_X\times\overline\Omega_Y)\times(\overline\Omega_X\times\overline\Omega_Y)\to\mathbb{R}$. Given such a function $\tilde\phi$, the standard double variable technique is to consider the function

\[
\phi:((x_1,y_1), (x_2,y_2)\mapsto u(x_1,y_1)-v(x_2,y_2)-\tilde\phi((x_1,y_1),(x_2,y_2))
\]

\noindent
and a point of maximum $((\overline x_1,\overline y_1),(\overline x_2,\overline y_2))$ of it. Such a maximum point gives the candidate points for the above explained procedure with

\begin{equation}
\label{eq:tests}
\begin{array}{ll}
\displaystyle
\varphi_1:(x_1,y_1)\mapsto\tilde\phi((x_1,y_1),(\overline x_2,\overline y_2)),\\
\displaystyle
\varphi_2:(x_2,y_2)\mapsto-\tilde\phi((\overline x_1,\overline y_1),(x_2,y_2)).
\end{array}
\end{equation}

As for the state-constraint optimal control problem, the main difficult is to be able to guarantee that, in the points $(\overline x_1,\overline y_1),(\overline x_2,\overline y_2)$, both equation inequalities hold (the last two lines of (\ref{eq:explanation})). This is made by a suitable use of the boundary conditions (\ref{eq:subsol}),(\ref{eq:supersol}), and by the use of a suitable penalizing term inside $\tilde\phi$ which avoids $(\overline x_1,\overline y_1),(\overline x_2,\overline y_2)$, or some components of them, to belong to the boundaries.

For this reason, in the following proof, we are mostly going to put in evidence the suitable penalizing term and the use of the boundary conditions for our particular case of problem (\ref{eq:subsol}), (\ref{eq:supersol}), and to refer to Bardi-Capuzzo Dolcetta \cite{barcap} (pages 279--280) for the "standard" part in order to obtain all the lines of (\ref{eq:explanation}) and the contradiction (\ref{eq:contradiction}). The point c) of the proof is more detailed.

}

{\it Proof of Theorem \ref{thm:uniqueness}.} 
\noindent
Let us assume the absurd hypothesis (\ref{eq:absurd}) and take $\delta>0$ and $(x_0,y_0)\in\overline\Omega_X\times\overline\Omega_Y$ such that the other following absurd hypothesis holds

\begin{equation}
\label{eq:absurd2}
u(x_0,y_0)-v(x_0,y_0)>m-\delta\ge\frac{m}{2}>0.
\end{equation}

\noindent
We only treat the boundary case $(x_0,y_0)\in\partial(\Omega_X\times\Omega_Y)$ and, if $(x_0,y_0)\in\partial\Omega_X\times\partial\Omega_Y$, $\psi_Y(x_0,y_0)<\psi_X(x_0,y_0)$. The other cases are similar or easier. \CorrB{For example, when $\psi_Y(x_0,y_0)=\psi_{XY}(x_0,y_0)=\psi_X(x_0,y_0)$, then we are in a sort of continuity case, all exit costs almost coincide around $(x_0,y_0)$, and the situation is similar to the one with continuous datum on the whole boundary.}

We take $\varepsilon>0$ and, for every one of the following cases, we consider a suitable ``double variable" function $\phi:(\overline\Omega_X\times\overline\Omega_Y)\times(\overline\Omega_X\times\overline\Omega_Y)\to\mathbb{R}$,
as here explained. In the following, $\zeta$ is a $C^1$ positive function on $\mathbb{R}^n\times\mathbb{R}^m$ with bounded gradient and such that $\zeta(x_0,y_0)=0$ and that $\zeta\to+\infty$ when $\|(x,y)\|\to+\infty$, and $\mu>0$ is a constant whose value will be fixed later \CorrB{(see Bardi-Capuzzo Dolcetta \cite{barcap}, page 54, for the use of this kind of functions in the comparison results in case of an unbounded domain)}.

a) i) $(x_0,y_0)\in\partial\Omega_X\times\Omega_Y$ and $v(x_0,y_0)<\psi_X(x_0,y_0)$; or ii) $(x_0,y_0)\in\Omega_X\times\partial\Omega_Y$ and $v(x_0,y_0)<\psi_Y(x_0,y_0)$; or iii) $(x_0,y_0)\in\partial\Omega_X\times\partial\Omega_Y$ and $\psi_Y(x_0,y_0)<v(x_0,y_0)<\psi_X(x_0,y_0)$; or iv) $(x_0,y_0)\in\partial\Omega_X\times\partial\Omega_Y$ and $\psi_Y(x_0,y_0)\neq v(x_0,y_0)<\psi_X(x_0,y_0)$:

\begin{equation}
\label{eq:doublevariable-a}
\begin{array}{l}
\displaystyle
\phi_a((x_1,y_1),(x_2,y_2))=u(x_1,y_1)-v(x_2,y_2)\\
\displaystyle
-\left\|\frac{x_1-x_2}{\varepsilon}-\eta_X(x_0)\right\|^2-\|x_2-x_0\|^2
-\left\|\frac{y_1-y_2}{\varepsilon}-\eta_Y(y_0)\right\|^2-\|y_2-y_0\|^2\\
\displaystyle
-\mu\zeta(x_1,y_1)-\mu\zeta(x_2,y_2).
\end{array}
\end{equation}

b) i) $(x_0,y_0)\in\partial\Omega_X\times\Omega_Y$ and $v(x_0,y_0)\ge\psi_X(x_0,y_0)$; or ii) $(x_0,y_0)\in\Omega_X\times\partial\Omega_Y$ and $v(x_0,y_0)\ge\psi_Y(x_0,y_0)$; or iii) $(x_0,y_0)\in\partial\Omega_X\times\partial\Omega_Y$ and $v(x_0,y_0)\ge\psi_X(x_0,y_0)$:

\begin{equation}
\label{eq:doublevariable-b}
\begin{array}{l}
\displaystyle
\phi_b((x_1,y_1),(x_2,y_2))=u(x_1,y_1)-v(x_2,y_2)\\
\displaystyle
-\left\|\frac{x_2-x_1}{\varepsilon}-\eta_X(x_0)\right\|^2-\|x_1-x_0\|^2
-\left\|\frac{y_2-y_1}{\varepsilon}-\eta_Y(y_0)\right\|^2-\|y_1-y_0\|^2\\
\displaystyle
-\mu\zeta(x_1,y_1)-\mu\zeta(x_2,y_2).
\end{array}
\end{equation}

c) $(x_0,y_0)\in\partial\Omega_X\times\partial\Omega_Y$ and $v(x_0,y_0)=\psi_Y(x_0,y_0)$:

\begin{equation}
\label{eq:doublevariable-c}
\begin{array}{ll}
\displaystyle
\phi_c((x_1,y_1),(x_2,y_2))=u(x_1,y_1)-v(x_2,y_2)\\
\displaystyle
-\left\|\frac{x_1-x_2}{\varepsilon}-\eta_X(x_0)\right\|^2-\|x_2-x_0\|^2
-\left\|\frac{y_2-y_1}{\varepsilon}-\eta_Y(y_0)\right\|^2-\|y_1-y_0\|^2\\
\displaystyle
-\mu\zeta(x_1,y_1)-\mu\zeta(x_2,y_2).
\end{array}
\end{equation}

Note the differences: from $\phi_a$ to $\phi_b$:  in all penalizing terms the role of indexes $1$ and $2$ are mutually exchanged; from $\phi_a$ to $\phi_c$: in the second penalizing terms the role of indexes $1$ and $2$ are mutually exchanged. This means that, when performing the usual double variable technique for comparison results, for suitable test functions $\varphi_1$ and $\varphi_2$ defined as in (\ref{eq:tests}), in the case a) we are going to detach  maxima for $u-\varphi_1$ (i.e. $(x_1,y_1)$) from the boundary; in the case b) we are going to detach minima for $v-\varphi_2$ (i.e. $(x_2,y_2)$) from the boundary; in the case c) we are going to detach the $x$-component of the maxima for $u-\varphi_1$ (i.e. $x_1$) and the $y$-component of the minima for $v-\varphi_2$ (i.e. $y_2$) from the boundary.  We briefly treat some of the above cases.

\CorrB{In the following $C>0$ is a suitable constant and $\omega$ is a suitable modulus of continuity, whose choices are independent from $\delta$ and $\varepsilon$ (depending on $\mu$ only). 
}

Case a). \CorrB{We use the absurd hypothesis (\ref{eq:absurd2}), the hypothesis on $\zeta$, the penalizing terms involving $\eta_X$ and $\eta_Y$, and the hypothesis (\ref{eq:cone}). Using standard estimates, we have that, at least for small $\delta$ and $\varepsilon$, $\phi_a$ has a maximum point in $(\overline\Omega_X\times\overline\Omega_Y)\times(\overline\Omega_X\times\overline\Omega_Y)$, let us say $((x_1^\varepsilon,y_1^\varepsilon),(x_2^\varepsilon,y_2^\varepsilon))$, and that

\[
\begin{array}{ll}
\displaystyle
\|x_1^\varepsilon-x_2^\varepsilon\|,\|y_1^\varepsilon-y_2^\varepsilon\|\le C\varepsilon,\\
\displaystyle
\|x_2^\varepsilon-x_0\|,\|y_2^\varepsilon-y_0\|\le\sqrt{\delta+\omega(\varepsilon)}\\
\displaystyle
\left\|\frac{x_1-x_2}{\varepsilon}-\eta_X(x_0)\right\|, \left\|\frac{y_1-y_2}{\varepsilon}-\eta_Y(y_0)\right\|\le\sqrt{\delta+\omega(\varepsilon)}\\
\displaystyle
(x_1^\varepsilon,y_1^\varepsilon)\in\Omega_X\times\Omega_Y.
\end{array}
\]
}

\noindent
i) For small $\varepsilon$ we have $(x_2^\varepsilon,y_2^\varepsilon)\in\overline\Omega_X\times\Omega_Y$ and, if $x_2^\varepsilon\in\partial\Omega_X$, then $v(x_2^\varepsilon,y_2^\varepsilon)<\psi_X(x_2^\varepsilon,y_2^\varepsilon)$. Hence both equation inequalities hold in (\ref{eq:subsol}) and in (\ref{eq:supersol}) for $u$ and $v$ respectively, when tested with test functions obtained from $\phi_a$ as in (\ref{eq:tests}). \CorrB{We then get the required contradiction, because, by the absurd hypothesis (\ref{eq:absurd2}) and the infinitesimal estimates here above, we get the first and the fourth lines of (\ref{eq:explanation}), just taking $\mu,\delta,\varepsilon$ sufficiently small.}

\noindent
The points ii), iii) and iv) are similarly treated.

b) As before, let $((x_1^\varepsilon,y_1^\varepsilon),(x_2^\varepsilon,y_2^\varepsilon))$ be a point of maximum for $\phi_b$. In this case, for at least small $\delta$ and $\varepsilon$, it is

\CorrB{
\[
\begin{array}{ll}
\displaystyle
\|x_1^\varepsilon-x_2^\varepsilon\|,\|y_1^\varepsilon-y_2^\varepsilon\|\le C\varepsilon,\\
\displaystyle
\|x_1^\varepsilon-x_0\|,\|y_1^\varepsilon-y_0\|\le\sqrt{\delta+\omega(\varepsilon)}\\
\displaystyle
\left\|\frac{x_2-x_1}{\varepsilon}-\eta_X(x_0)\right\|, \left\|\frac{y_2-y_1}{\varepsilon}-\eta_Y(y_0)\right\|\le\sqrt{\delta+\omega(\varepsilon)}\\
\displaystyle
(x_2^\varepsilon,y_2^\varepsilon)\in\Omega_X\times\Omega_Y.
\end{array}
\]
}

i) The hypothesis $v(x_0,y_0)\ge\psi_X(x_0,y_0)$ and the absurd hypothesis (\ref{eq:absurd2}) imply $u(x_0,y_0)>\psi_X(x_0,y_0)$ and so $u(x_1^\varepsilon,y_1^\varepsilon)>\psi_X(x_1^\varepsilon,y_1^\varepsilon)$, for small $\delta$ and $\varepsilon$. Hence both equation inequalities hold in (\ref{eq:subsol}) and in (\ref{eq:supersol}) for $u$ and $v$ respectively, when tested with suitable test functions obtained from $\phi_b$ as in (\ref{eq:tests}). We then get the contradiction as before.

\noindent
The points ii) and iii) are similarly treated.

c) Again, let $((x_1^\varepsilon,y_1^\varepsilon),(x_2^\varepsilon,y_2^\varepsilon))$ be a point of maximum for $\phi_c$. For at least small $\delta$ and $\varepsilon$, it is

\CorrB{
\[
\begin{array}{ll}
\displaystyle
\|x_1^\varepsilon-x_2^\varepsilon\|,\|y_1^\varepsilon-y_2^\varepsilon\|\le C\varepsilon,\\
\displaystyle
\|x_2^\varepsilon-x_0\|,\|y_1^\varepsilon-y_0\|\le\sqrt{\delta+\omega(\varepsilon)}\\
\displaystyle
\left\|\frac{x_1-x_2}{\varepsilon}-\eta_X(x_0)\right\|, \left\|\frac{y_2-y_1}{\varepsilon}-\eta_Y(y_0)\right\|\le\sqrt{\delta+\omega(\varepsilon)}\\
\displaystyle
x_1^\varepsilon\in\Omega_X,\ y_2^\varepsilon\in\Omega_Y.
\end{array}
\]
}

\noindent
The hypothesis $v(x_0,y_0)=\psi_Y(x_0,y_0)$ and absurd hypothesis (\ref{eq:absurd2}) imply that, for small $\delta$ and $\varepsilon$, $u(x_1^\varepsilon,y_1^\varepsilon)>\psi_Y(x_1^\varepsilon,y_1^\varepsilon)$ and $v(x_2^\varepsilon,y_2^\varepsilon)<\psi_X(x_2^\varepsilon,y_2^\varepsilon)$.

For a seek of completeness, we show here the calculation for this case. In particular, let us note that in this case we are not detaching from the boundary both points of maximum and of minimum, which is in general not possible, but instead we are detaching the first $n$ components of the point of maximum and the second $m$ components of the point of minimum. And this is possible because the domain $\Omega=\Omega_X\times\Omega_Y\subset\mathbb{R}^n\times\mathbb{R}^m$ is indeed a cartesian product. By definition of $\phi_c$ we have

\[
\phi_c((x_0,y_0),(x_0,y_0))=u(x_0,y_0)-v(x_0,y_0)-\|\eta_X(x_0)\|-\|\eta_Y(y_0)\|\ge m-\delta-2M,
\]

\noindent
where $M$ is a bound for $\|\eta_X\|$ and for $\|\eta_Y\|$ (and also for $|u|$ and $|v|$, in subsequent calculations). By the coercivity of $\zeta$ and the boundedness of $u$ and $v$ we get that $\phi_c$ reaches its maximum in a point $((x_1^\varepsilon,y_1^\varepsilon),(x_2^\varepsilon,y_2^\varepsilon))$ and that there exists two compact subsets (depending on $\mu$) $K^\mu_X\subseteq\overline\Omega_X$, $K^\mu_Y\subseteq\overline\Omega_Y$ such that $((x_1^\varepsilon,y_1^\varepsilon),(x_2^\varepsilon,y_2^\varepsilon))\in(K^\mu_X\times K_Y)\times(K^\mu_X\times K_Y)$, for every  $\varepsilon$ sufficiently small. We can also suppose $x_0+\varepsilon\eta_X(x_0)\in K^\mu_X$, $y_0+\varepsilon\eta_Y(y_0)\in K^\mu_Y$ for all $\varepsilon>0$ sufficiently small. Let $\omega^\mu$ be a modulus of continuity for both $u$ and $v$ and for $\zeta$ in $(K^\mu_X\cap\overline\Omega_X)\times(\times K^\mu_Y\cap\overline\Omega_Y)$. We have

\begin{equation}
\label{eq:ineq1}
\begin{array}{l}
\displaystyle
\phi_c((x_1^\varepsilon,y_1^\varepsilon),(x_2^\varepsilon,y_2^\varepsilon))\ge
\phi_c((x_0+\varepsilon\eta_X(x_0),y_0),(x_0,y_0+\varepsilon\eta_Y(y_0))\\
\displaystyle
=u(x_0+\varepsilon\eta_X(x_0),y_0)-v(x_0,y_0+\varepsilon\eta_Y(y_0))\\
\displaystyle
-\zeta(x_0+\varepsilon\eta_X(x_0),y_0)-\zeta(x_0,y_0+\varepsilon\eta_Y(y_0))
\ge m-\delta-4\omega^\mu(C\varepsilon),
\end{array}
\end{equation}

\noindent
where $C>0$ is a suitable constant independent from $\varepsilon$. Now, we have the inequalities, for $\delta$ and $\varepsilon$ small,

\begin{equation}
\label{eq:ineq2}
\begin{array}{l}
\displaystyle
u(x_1^\varepsilon,y_1^\varepsilon)-v(x_2^\varepsilon,y_2^\varepsilon)\le2M,\\
\displaystyle
0<\frac{m}{2}<u(x_1^\varepsilon,y_1^\varepsilon)-v(x_2^\varepsilon,y_2^\varepsilon)\le m+\omega^\mu(\|(x_1^\varepsilon,y_1^\varepsilon)-(x_2^\varepsilon,y_2^\varepsilon)\|).
\end{array}
\end{equation}

\noindent
From the definition of $\phi_c$ and from (\ref{eq:ineq1}), we get

\begin{equation}
\label{eq:ineq3}
\begin{array}{l}
\displaystyle
\left\|\frac{x_1^\varepsilon-x_2^\varepsilon}{\varepsilon}-\eta_X(x_0)\right\|^2+\|x_2^\varepsilon-x_0\|^2
+\left\|\frac{y_2^\varepsilon-y_1^\varepsilon}{\varepsilon}-\eta_Y(y_0)\right\|^2+\|y_1^\varepsilon-y_0\|^2\\
\displaystyle
\le\delta+\omega^\mu(\|(x_1^\varepsilon,y_1^\varepsilon)-(x_2^\varepsilon,y_2^\varepsilon)\|)+4\omega^\mu(C\varepsilon).
\end{array}
\end{equation}

\noindent
By the boundedness of $\omega^\mu$, when its argument is the distance of points in $K_X^\mu\times K_Y^\mu$, and the boundedness of $\eta_X,\eta_Y$, from (\ref{eq:ineq3}) we get (for another constant independent from $\varepsilon$, and still denoted by $C$)

\begin{equation}
\label{eq:ineq3.1}
\|x_1^\varepsilon-x_2^\varepsilon\|+\|y_1^\varepsilon-y_2^\varepsilon\|\le C\varepsilon,
\end{equation}

\noindent
which, again by (\ref{eq:ineq3}) and for another $C>0$ independent from $\varepsilon$, gives

\begin{equation}
\label{eq:ineq4}
\begin{array}{l}
\displaystyle
\left\|\frac{x_1^\varepsilon-x_2^\varepsilon}{\varepsilon}-\eta_X(x_0)\right\|^2+\|x_2^\varepsilon-x_0\|^2
+\left\|\frac{y_2^\varepsilon-y_1^\varepsilon}{\varepsilon}-\eta_Y(y_0)\right\|^2+\|y_1^\varepsilon-y_0\|^2\\
\displaystyle
\le \delta+5\omega^\mu(C\varepsilon).
\end{array}
\end{equation}
By the inward cone hypothesis (\ref{eq:cone}), from (\ref{eq:ineq4}), in a standard way (Bardi-Capuzzo Dolcetta \cite{barcap}, page 280), we get, for sufficiently small $\delta$ and $\varepsilon$,
\begin{equation}
\label{eq:interior}
x_1^\varepsilon\in\Omega_X,\ y_2^\varepsilon\in\Omega_Y.
\end{equation}

Now, since we are in the case $v(x_0,y_0)=\psi_Y(x_0,y_0)<\psi_X(x_0,y_0)$  and since $u(x_0,y_0)>v(x_0,y_0)$, by (\ref{eq:ineq3.1})--(\ref{eq:ineq4}) we can suppose that 

\[
 u(x_1^\varepsilon,y_1^\varepsilon)>\psi_Y(x_1^\varepsilon,y_1^\varepsilon),\ v(x_2^\varepsilon,y^2_\varepsilon)<\psi_X(x_2^\varepsilon,y_2^\varepsilon).
\]
\noindent
This, together with (\ref{eq:interior}) and the definition of $((x_1^\varepsilon,y_1^\varepsilon),(x_2^\varepsilon,y_2^\varepsilon))$, implies that both equation inequalities hold, in $(x_1^\varepsilon,y_1^\varepsilon)$ for $u$ as in (\ref{eq:subsol}) and in $(x_2^\varepsilon,y_2^\varepsilon)$ for $v$ as in (\ref{eq:supersol}) respectively, when we take as test functions

\[
\begin{array}{l}
\displaystyle
\varphi_1(x,y)=\left\|\frac{x-x_2^\varepsilon}{\varepsilon}-\eta_X(x_0)\right\|^2+\left\|\frac{y_2^\varepsilon-y}{\varepsilon}-\eta_Y(y_0)\right\|^2+\|y-y_0\|^2+\mu\zeta(x,y)\\
\displaystyle
\varphi_2(x,y)=-\left\|\frac{x_1^\varepsilon-x}{\varepsilon}-\eta_X(x_0)\right\|^2-\|x-x_0\|^2
-\left\|\frac{y-y_1^\varepsilon}{\varepsilon}-\eta_Y(y_0)\right\|^2-\mu\zeta(x,y),\\
\end{array}
\]

\noindent
respectively. We have

\[
\begin{array}{l}
\displaystyle
\nabla\varphi_1(x^\varepsilon_1,y^\varepsilon_1)\\
\displaystyle
=\mu\nabla\zeta(x_1^\varepsilon,y_1^\varepsilon)+2(0,y_1^\varepsilon-y_0)
+\frac{2}{\varepsilon}\left(\frac{x^\varepsilon_1-x^\varepsilon_2}{\varepsilon}-\eta_X(x_0),\frac{y^\varepsilon_1-y^\varepsilon_2}{\varepsilon}-\eta_Y(y_0)\right),\\
\displaystyle
\nabla\varphi_2(x^\varepsilon_2,y^\varepsilon_2)\\
\displaystyle
=-\mu\nabla\zeta(x_2^\varepsilon,y_2^\varepsilon)-2(x_2^\varepsilon-x_0,0)
+\frac{2}{\varepsilon}\left(\frac{x^\varepsilon_1-x^\varepsilon_2}{\varepsilon}-\eta_X(x_0),\frac{y^\varepsilon_1-y^\varepsilon_2}{\varepsilon}-\eta_Y(y_0)\right),\\
\end{array}
\]

\noindent
and then, if $\mu$ is sufficiently small, we can conclude in the standard way getting the conclusion by contradiction to (\ref{eq:absurd}). 
\cvd

\begin{Remark}
\label{rmrk:Nash}
As already remarked, the exit cost $\psi_{XY}$, for simultaneous exit of $X$ and $Y$, does not play any role in the formulation of the Dirichlet problem (\ref{eq:isaacs}). Indeed, it can never happen that the simultaneous exit cost $\psi_{XY}$ is a ``good" choice for both players (i.e. an equilibrium) without being already equal to $\psi_X$ or to $\psi_Y$ or to $V_{int}^\delta$ for some $\delta>0$  where the latter is defined as the lower value function restricted to controls $\beta$ and to non-anticipating strategies $\gamma$ which make $Y$ and $X$ remain inside $\overline\Omega_Y$ and $\overline\Omega_X$ for times in $[0,\delta]$, respectively.  For instance, let us suppose that, in a point $(x,y)\in\partial\Omega_X\times\partial\Omega_Y$, we have

\[
\psi_Y<\psi_{XY}<\min\{\psi_X,V_{int}^\delta\}.
\]

\noindent
Then, player II (the maximizing one) has certainly no interest in exit, and so the ``really paid cost" is $\psi_X$ or $V_{int}^\delta$. A similar conclusion holds for the case

\[
\max\{\psi_Y,V_{int}^\delta\}<\psi_{XY}<\psi_X.
\] 

\noindent
In the case that $\psi_{XY}=V_{int}^\delta$ is a "good choice" for both players, then dynamic programming leads to the Isaacs equation and so the exit cost $\psi_{XY}$ does not really influence the problem.

By the way, even in a strategic static minmax game where two players may independently choose to "stay" or to "exit" and the first player wants to minimize, if the utility $u(exit,exit)$ stays between the utilities $u(stay,exit)\le u(exit,stay)$, then the choice $(exit,exit)$ is never a Nash equilibrium, whichever $u(stay,stay)$ is.

\end{Remark}

\begin{Remark}
Since the dynamics are decoupled, in order to have the classical Isaacs' condition for the existence of a value of the game (see for example Bardi-Capuzzo Dolcetta \cite{barcap}), we only need some further hypotheses on the running cost $\ell$. The simplest one is that it is also decoupled with respect to controls \CorrR{(i.e. $\ell(x,y,a,b)=\ell_1(x,y)+\ell_2(a)+\ell_3(b)$) }. In this case the two Hamiltonians $UH$ and $LH$ are the same and hence, by uniqueness of the corresponding Dirichlet problems (\ref{eq:isaacs}), (\ref{eq:isaacs2}), $\overline V=\underline V$. 
\end{Remark}

\section{On constrained non-anticipating strategies}
\label{sec:constrained_strategies}

\CorrR{We give a possible construction of state-constraint non-anticipating tuning as well as non-anticipating strategies satisfying Assumption 2. We follow Soner's \cite{son} construction of state-constraint controls (see also Bardi-Capuzzo Dolcetta \cite{barcap}, pages 272--274), modifying it in a non-anticipating way, in order to adapt such a construction to our purposes (see Remark \ref{rmrk:soner}).}

We are now considering only point I) of Assumption 2, point II) being similar. Note that point I) is concerning with non-anticipating strategies $\gamma\in\Gamma$ for player $X$. With respect to (\ref{eq:hypotheses}) (and to (\ref{eq:systems})), we are going to relax a little bit the hypotheses of decoupled dynamics and we are going to consider the following hypothesis:  the dynamics $f$, the one for player $X$, is affine with respect to the controls (coherently with Bettiol-Cardialaguet-Quincampoix \cite{betcarqui}) and ``weakly'' depends on the control of $Y$. \CorrR{This means that (here, $\tilde f$ stays for the dynamics of the first player $X$)

	\begin{equation}\label{nuovaf}
	\begin{array}{l}
\displaystyle
\tilde f:\mathbb{R}^n\times A\times B\to\mathbb{R}^n,\ f(x, a , b)=f(x)+a+Db,\\
\displaystyle
A\subset\mathbb{R}^n,\ DB\subseteq A\ \mbox{are compact sets}
	\end{array}
	\end{equation}
where $D \in \mathbb{R}^{n' \times m'}$ is a fixed constant matrix, $B$ is as in (\ref{eq:hypotheses}), and $f:\mathbb{R}^n\to\mathbb{R}^n$ is bounded and Lipschitz continuous.  Note that, the dynamics $\tilde f$ then satisfies similar regularity hypotheses as in (\ref{eq:hypotheses}), in particular, for some $M,L>0$ and for all $x_1,x_2\in\mathbb{R}^n,a\in A,b\in B$:

\[
\|\tilde f(x,a,b)\|\le M,\ \ \|\tilde f(x_1,a,b)-\tilde f(x_2,a,b)\|\le L\|x_1-x_2\|.
\]

 \noindent
The system for the trajectories of the first player is

\[
\left\{
\begin{array}{ll}
\displaystyle
X'(t)=f(X(t))+\alpha(t)+D\beta(t)\\
\displaystyle
X(0)=x\in\overline\Omega_X
\end{array}
\right.
\]

\begin{Remark}
\label{rmrk:weakly}
We point out that here and in what follows, we are assuming such a weak decoupling for the dynamics of the player $X$ only. The dynamics $g$ of the player $Y$ will be still considered decoupled, i.e. depending on the pair $(y,b)$ only, as in (\ref{eq:hypotheses}). In this situation, we are going to prove that all the hypotheses i)--vi) of item I) only, of Assumption 2, are satisfied. The fact that the dynamics of the second player is completed decoupled, i. e. it does not depend on the control of the first player, enters in what follows because $\tilde\beta$ in (\ref{betatilde}) is constructed independently on the behavior and on the controls of the first player.

On the other hand, if we assume that the dynamics $g$ is weakly decoupled (similarly to (\ref{nuovaf}), changing the role of $a$ and $b$), and we maintain the decoupled feature of $f$ as in (\ref{eq:hypotheses}), then we can prove that item II) of Assumption 2 is satisfied. 

Also note that, the validity of item I) (respectively, item II)) permits to conclude that the lower value function $\underline V$ (respectively, the upper value function $\overline V$) is continuous. However, all the proofs of the results in Sections \ref{sec:isaacs} and \ref{sec:uniqueness} (as well as the simultaneous validity of both items I) and II))
hold in the case of decoupling of both $f$ and $g$ as in (\ref{eq:hypotheses}), that is when $D=0$ in (\ref{nuovaf}). The extension of such results to the weak decoupled case may be the subject of future works.

\CorrRR{However, we point out that, in order to get the estimate on the costs, we are going to also assume a decoupled feature of the running cost with respect to the controls (see (\ref{sumrunningcosts}))}.

\end{Remark}

Assuming (\ref{nuovaf}), we need a modification of the controllability Assumption 1.}

{\it Assumption 3}. Similarly as in Assumption 1, we assume here that, for every $x\in\partial\Omega_X$, there exist two constant controls $a_1,a_2\in A$ such that  $\tilde f(x,a_1, b)$ is {\it strictly entering in $\Omega_X$} and $\tilde f(x,a_2, b)$ is {\it strictly entering in $\mathbb{R}^m\setminus\overline\Omega_X$} $\forall b \in B$.\\

\CorrRR{Assuming $C^2$-regularity of $\partial\Omega_X$ (\ref{eq:hypotheses}), by Assumption 3, by the weak decoupling (\ref{nuovaf}), and by the Lipschitz continuity of $f$, for every compact $K$, there exist $\zeta>0$ and $r>0$ and, for any $\overline x \in K\cap\partial\Omega_X$, there exists $a(\overline x)\in A$ such that, for every $x\in B(\overline x,r)\cap\overline\Omega_X$,
\begin{equation}\label{fscalarexi}
\inf_{b\in B}(f(x)+a(\overline x)+Db)\cdot \xi(\overline x)>0,\ \ \ (f(x)+a(\overline x))\cdot\xi(\overline x)\ge\zeta
\end{equation}
}

\noindent
where $\xi(x)$ is the inward normal unit vector to $\Omega_X$ at $x\in\partial\Omega_X$. In what follows, in view of possible future applications to thermostatically switching systems, we assume that the boundaries of $\Omega_X$ and $\Omega_Y$ are hyperplanes passing through the origin, and that $\Omega_X$ and $\Omega_Y$ are just one of the two semi-space defined by the hyperplane. In this way, the unit vector inward normal $\xi$ is constant on $\partial\Omega_X$ (as well as on $\partial\Omega_Y$, in the sequel denoted by the same letter $\xi$). However, everything done here can be easily generalized to the case in which the boundaries of $\Omega_X$ and $\Omega_Y$ are finite intersections of hyperplanes not necessarily passing through the origin (see for example what done in Bagagiolo-Bardi \cite{bagbar}). Moreover it can be extended to more general regular domains.

\CorrR{Now, we are going to prove that, in the situation described above,  all the points i)--vi) of item I) of Assumption 2 are satisfied.} \CorrRR{In what follows, we will indicate by $K$ a generic compact set of the form $K=K^X\times K^Y$, where $K^X\subseteq\overline\Omega_X$ and $K^Y\subseteq\overline\Omega_Y$ are compact.}

i), iii) and v). Take $T>0$ , and take $T\ge t_Y^*>0$ to be fixed later on and $y_1,y_2\in \CorrRR{K^Y}\subseteq\overline\Omega_Y$ compact, and define

\begin{equation}\label{defepsilonY}
\varepsilon_Y=\sup_{\beta\in\mathcal{B}}\left(\sup_{0\le t\le\min\{\tau_Y(y_1,\beta),t_Y^*\}} \left(-\xi\cdot  Y(t;y_2,\beta)\right)^+\right),
\end{equation}
where $(r)^+=\max(r,0)$ is the positive part, and $\xi$ is the unit internal normal to $\Omega_Y$. Note that, since $\Omega_Y$ is a semi-space, the quantity inside the suprema over $\beta$ is just the maximal distance from $\overline\Omega_Y$ reached by the trajectory starting from $y_2$ with control $\beta$, before that the trajectory starting from $y_1$ with the same control  $\beta$ exits from $\overline\Omega_Y$, or the time $t^*_Y$ is reached. \CorrR{However, the presence of the supremum over $\beta$ makes $\varepsilon_Y$ independent from $\beta$, and this is the essential feature for the fact that next formula (\ref{betatilde}) defines a non-anticipating tuning, as we are going to explain in the comments after (\ref{betatilde}) and in Remark \ref{rmrk:soner}. }

Note that we have the estimate (with $C$ depending only on $T$, $K$, and $t_Y^*$)
\begin{equation}
\label{eq:epsilonestimate}
0\le\varepsilon_Y\le C(t_Y^*)\|y_1-y_2\|.
\end{equation}

\noindent
Indeed, for every $\beta$ and for every $0\le t\le\min\{\tau_Y(y_1,\beta),t_Y^*\}$, we have $Y(t;y_1,\beta)\in\overline\Omega_Y$, i.e. $\xi \cdot Y(t;y_1,\beta)\geq 0$. Then there exists a constant $C(t_Y^*)$ such that 
\begin{equation}
\begin{array}{l}
\displaystyle
-\xi \cdot Y(t;y_2,\beta)\leq -\xi \cdot Y(t;y_2,\beta) - (-\xi\cdot Y(t;y_1,\beta))\\
\displaystyle
\leq  \|Y(t;y_2,\beta)-Y(t;y_1,\beta)\|\leq C(t_Y^*)\|y_1-y_2\|,
\end{array}
\end{equation}
and hence (\ref{eq:epsilonestimate}) holds.

\noindent
Now, take $\beta\in\mathcal{B}$, let $t_{0_Y}\ge0$ be the first time the trajectory $Y(\cdot):=Y(\cdot;y_2,\beta)$ hits the boundary $\partial\Omega_Y$ and let $b_0\in B$ be such that $g(Y(t_{0_Y}),b_0)$ strictly enters in $\Omega_Y$ at $Y(t_{0_Y})$: the one given by Assumption 1. Now, let us take $k_Y>0$ and define the measurable control $\tilde\beta\in\mathcal{B}$ as

\begin{equation}
\label{betatilde}
\tilde\beta(t)=\left\{
\begin{array}{ll}
\displaystyle
\beta(t)&\mbox{if } 0\le t\le\min\{t_{0_Y},t_Y^*\},\\
\displaystyle
b_0&\mbox{if } \min\{t_{0_Y},t_Y^*\}\le t\le\min\{t_{0_Y},t_Y^*\}+k_Y\varepsilon_Y,\\
\displaystyle
\beta(t-k_Y\varepsilon_Y)&\mbox{if } t\ge\min\{t_{0_Y},t_Y^*\}+k_Y\varepsilon_Y 
\end{array}
\right.
\end{equation}

\CorrR{By our definition of $\varepsilon_Y$ (\ref{defepsilonY}) which is independent on $\beta$, due to the presence of the supremum over the controls, the construction in (\ref{betatilde}) is a non-anticipating tuning (see Definition \ref{def:nonant}, see also Remark \ref{rmrk:soner} for other comments). This means that, whenever $\beta_1=\beta_2$ in $[0,t]$ a.e., then also $\tilde\beta_1=\tilde\beta_2$ in $[0,t]$ a.e., in other words  i) holds}. Indeed, in such a case, in the time interval $[0,t]$, $\beta_1$ and $\beta_2$ generate the same trajectory $Y$ starting from $y_2$. If, in the time interval $[0,t]$ the trajectory $Y$ does not hit the boundary, then, by definition (25), $\tilde\beta_1=\beta_1=\beta_2=\tilde\beta_2$ in $[0,t]$. If instead the trajectory $Y$ hits the boundary at $t_{0_Y}\le t$, then in the interval $[0,t_{0_Y}[$ we still have the equality $\tilde\beta_1=\tilde\beta_2$, and in the time interval
$[t_{0_Y},t_{0_Y}+k_Y\varepsilon_Y]\cap[0,t]$ \CorrR{(whose length by (\ref{defepsilonY}) is independent from $\beta_1$ and $\beta_2$)}, we have $\tilde\beta_1=b_0=\tilde\beta_2$, with $b_0$ given by Assumption 1 in $Y(t_{0_Y})$, which is the same for both controls. Finally, in the (possibly empty) time interval $[0,t]\setminus[t_{0_Y},t_{0_Y}+k_Y\varepsilon_Y]$ we have $\tilde\beta_1(s)=\beta_1(s-k_Y\varepsilon_Y)=\beta_2(s-k_Y\varepsilon_Y)=\tilde\beta_2(s)$.

 Now we want to suitably choose $t_Y^*$ and $k_Y$ such that, for every $t\in[0,t_Y^*]$, it is $Y(t;y_2,\tilde\beta)\in\overline\Omega_Y$, at least for $t\le\tau_Y(y_1,\beta)$. This can be done (independently on $\beta$) just following Bardi-Capuzzo Dolcetta \cite{barcap} page 273, with $\varepsilon$ given by our $\varepsilon_Y$ (a more detailed construction is given for the similar question in the next point). Repeating the construction for every needed time-interval $[nt_Y^*,(n+1)t_Y^*]$, in order to cover the interval $[0,\tilde\tau]$,  we get iii). Finally, using (\ref{eq:epsilonestimate}), we also get v). 

\CorrR{
\begin{Remark}
\label{rmrk:soner}
Observe that (\ref{defepsilonY}) is different from the one defined by Soner \cite{son} (see also Bardi-Capuzzo Dolcetta \cite{barcap}, page 273)  since here we are building non-anticipating tuning as well as non-anticipating strategies, and this feature is guaranteed by the supremum over $\beta$ in definition (\ref{defepsilonY}). Indeed, in our framework and notations, if we just follow what done in \cite{son}, we would have, for every control $\beta$,

\[
\varepsilon(\beta)=\sup_{0\le t\le\min\{\tau_Y(y_1,\beta),t_Y^*\}} \left(-\xi\cdot  Y(t;y_2,\beta)\right)^+.
\]

\noindent
This means that, even if $\beta_1=\beta_2$ in $[0,t]$, then they may generate different values of 
$\varepsilon(\beta_1)$ and $\varepsilon(\beta_2)$ because they, and their corresponding trajectories, may differ after the time $t$. In particular, it may happen that the trajectories hit the boundary at the same instant $0\le t_{0_{Y}}<\min\{\tau_Y(y_1,\beta_1),\tau_Y(y_1,\beta_2),t^*_Y,t\}$ and that $[t_{0_Y},t_{0_Y}+k_Y\varepsilon(\beta_1)]\subset[t_{0_Y},t_{0_Y}+k_Y\varepsilon(\beta_2)]\subset[0,t]$. Applying (\ref{betatilde}), we would get two different behavior of $\tilde\beta_1$ and $\tilde\beta_2$ in $[0,t]$, that is an anticipating construction: (\ref{betatilde}) is not more a non-anticipating tuning (see Definition \ref{def:nonant}). Considering instead the supremum over all controls $\beta$, as we do in (\ref{defepsilonY}), makes us to avoid this behaviour because the length of $\varepsilon_Y$ does not depend by the single control $\beta$. 

Of course, as already said in the Introduction, in \cite{son} and \cite{barcap}, the non-anticipating structure is not taken into consideration because the reference problem is an optimal control problem.
Also observe that the construction (\ref{betatilde}) is exactly the same as in \cite{son} and \cite{barcap}. But here, the non-anticipating feature is given by the definition of $\varepsilon_Y$, which is independent on the single control. 

Finally observe that we may have $\varepsilon_Y=0$.  By its very definition (22) (do not consider here $t^*_Y$), $\varepsilon_Y=0$ means that, whatever the control $\beta$ is, the trajectory $Y(\cdot;y_2,\beta)$ does not exit from $\overline\Omega_Y$ before the trajectory $Y(\cdot;y_1,\beta)$ exits from $\overline\Omega_Y$. And this is exactly what we need in our proof of continuity of the value function (i. e.  requirement iii) of Assumption 2). Hence in this case, we do not need to modify the control: in (\ref{betatilde}), for every control $\beta$, it is $\tilde\beta=\beta$, coherently with the fact that the time interval of length $k_Y\varepsilon_Y$ (where we should make the modification of $\beta$) is just a point (not a true interval).

\end{Remark}
}

ii) and iv). Take $T>0$, and take $T\ge t_X^*$ to be fixed later on, $x_1, x_2 \in \CorrRR{K^X}\subseteq\overline\Omega_X$ compact and define, similarly as for $\varepsilon_Y$,

\begin{equation}\label{defepsilonX}
\varepsilon_X=\sup_{\gamma\in \Gamma}\sup_{\beta\in{\cal B}}\left(\sup_{0\le t\le\min\{\tau_X(x_2, \gamma[\CorrRR{\tilde\beta}], \CorrRR{\tilde\beta}),t_X^*\}}\left(-\xi\cdot  X(t;x_1,\gamma[\CorrRR{\tilde\beta}], \beta)\right)^+\right),
\end{equation}

where $\xi$ is the unit internal normal to $\Omega_X$, \CorrRR{and $\tilde\beta$ is defined as in (\ref{betatilde}), with respect to the previously fixed $y_1,y_2\in K^Y$.}  In (\ref{defepsilonX}) we use both the supremum over $\beta$ and over $\gamma$ in order to build non-anticipating strategies.
Note that now, in the notations of the trajectory $X$ and of the exit time $\tau_X$ we are taking account that the dynamics $f$ is only weakly decoupled (\ref{nuovaf}). \CorrRR{In particular, in this case, the estimates ii) and iv) of Assumption 2 must be replaced by (see also point 10) of the proof of Proposition \ref{prop:continuity})

\begin{equation}
\label{eq:newassumption}
\begin{array}{ll}
\displaystyle
ii')\ \tau_X(x_1,\tilde\gamma[\tilde\beta],\beta)\ge\tau_X(x_2,\gamma[\tilde\beta],\tilde\beta),\\
\displaystyle
iv')\ \|X(\tilde\tau;x_1,\tilde\gamma[\tilde\beta],\beta)-X(\tilde\tau;x_2,\gamma[\tilde\beta],\tilde\beta)\|\le{\cal O}_{T,K}(\|x_1-x_2\|+\|y_1-y_2\|)
\end{array}
\end{equation}

We refer to Remark \ref{rmrk:gluing_decoupling} for comments on the simpler case of strongly decoupled dynamics $f$.

}
As before, inside the suprema in (\ref{defepsilonX}), the scalar product is the distance from the semi-space $\overline\Omega_X$. \CorrRR{In this case we have (with $C$ depending only on $T$, $K$, $t_X^*$ and $t_Y^*$)
\begin{equation}\label{epsilonXestimate}
0\le\varepsilon_X\le C(\|x_1-x_2\|+\|y_1-y_2\|).
\end{equation}

\noindent
Indeed, for $t$ as in (\ref{defepsilonX}) and using (\ref{eq:epsilonestimate}),

\[
\begin{array}{ll}
\displaystyle
-\xi\cdot X(t;x_1,\gamma[\tilde\beta],\beta)\le\\
\displaystyle
-\xi\cdot X(t;x_1,\gamma[\tilde\beta],\beta)-(-\xi\cdot X(t;x_2,\gamma[\tilde\beta],\tilde\beta))\le\\
\displaystyle
\|x_1-x_2\|+L\int_0^t\|X(s;x_1,\gamma[\tilde\beta],\beta)-X(s;x_2,\gamma[\tilde\beta],\tilde\beta)\|ds+\int_0^tD(\beta(s)-\tilde\beta(s))ds\le\\
\displaystyle
\|x_1-x_2\|+L\int_0^t\|X(s;x_1,\gamma[\tilde\beta],\beta)-X(s;x_2,\gamma[\tilde\beta],\tilde\beta)ds+C(t_Y^*)\|y_1-y_2\|
\end{array}
\]

\noindent
and we conclude by the Gronwall estimate. Note that, in the estimate of the integral of $D(\beta-\tilde\beta)$, we have used the equality

\begin{equation}
\label{eq:betaestimate}
\begin{array}{ll}
\displaystyle
\int_0^tD(\beta(s)-\tilde\beta(s))ds=\\
\displaystyle
\int_{t_{0_Y}}^{t_{0_Y}+k_Y\varepsilon_Y}D(\beta(s)-b_0)ds+
\int_{t_{0_Y}+k_Y\varepsilon_Y}^tD(\beta(s)-\beta(s-k_Y\varepsilon_Y))ds=\\
\displaystyle
\int_{t_{0_Y}}^{t_{0_Y}+k_Y\varepsilon_Y}D(\beta(s)-b_0)ds-
\int_{t_{0_Y}}^{t_{0_Y}+k_Y\varepsilon_Y}D\beta(s)ds+
\int_{t-k_Y\varepsilon_Y}^tD\beta(s)ds
\end{array}
\end{equation}

}
Now, take $\gamma \in \Gamma$, and we want to construct  the strategy $\tilde\gamma$ acting on \CorrRR{$\tilde\beta$ (see (\ref{eq:newassumption}))} . Take $\beta\in{\cal B}$ and let $t_{0_X}\ge0$ be the first time the trajectory $X(\cdot;x_1,\gamma[\CorrRR{\tilde\beta}], \beta)$ hits the boundary $\partial\Omega_X$ and let $a_0\in A$ be such that $f$ strictly enters in $\Omega_X$ as in Assumption 3. Now, let us take $k_X>0$ and define \CorrRR{$\tilde{\tilde\gamma}:{\mathcal B}\to{\mathcal A}$} as

\begin{equation}
\label{gammatilde}
\begin{array}{l}
\displaystyle
\CorrRR{\tilde{\tilde{\gamma}}}[\beta](t)=\left\{
\begin{array}{ll}
\displaystyle
\gamma[\CorrRR{\tilde\beta}](t)&\mbox{if } 0\le t\le\min\{t_{0_X},t_X^*\},\\
\displaystyle
a_0&\mbox{if } \min\{t_{0_X},t_X^*\}\le t\le\min\{t_{0_X},t_X^*\}+k_X\varepsilon_X,\\
\displaystyle
\gamma[\CorrRR{\tilde\beta}](t-k_X\varepsilon_X)&\mbox{if } t\ge\min\{t_{0_X},t_X^*\}+k_X\varepsilon_X 
\end{array}
\right.
\end{array}
\end{equation}

\noindent
Being $\gamma$ \CorrRR{ and the tuning $\beta\mapsto\tilde\beta$} non-anticipating, by the definition of $\varepsilon_X$ we get that $\CorrRR{\tilde{\tilde \gamma}}$ is also non-anticipating. \CorrRR{We use the following notations, which also give notational coherence with (\ref{defepsilonX}) and (\ref{eq:newassumption}),

\[
\tilde\gamma[\tilde\beta](t)=\tilde{\tilde\gamma}[\beta](t),\ \ X(\cdot)=X(\cdot;x_1,\gamma[\tilde\beta],\beta),\ \ 
\tilde X(\cdot)=X(\cdot;x_1,\tilde\gamma[\tilde\beta],\beta).
\]

} 

\noindent
As before we want to suitably choose $t_X^*$ and $k_X$ such that, for every $t \in [0, t_X^*]$, it is $\tilde X(t)\in\overline\Omega_X$. Of course we are interested in the case $t\le\tau_X(x_2,\gamma[\CorrRR{\tilde\beta}],\CorrRR{\tilde\beta})$. We then prove that $$\xi \cdot \tilde X(t)\geq 0.$$
Again, we follow \cite{barcap} but, due to presence of both controls (\ref{nuovaf}), we now show some explicit calculations. 
Note that if $\min\{t_{0_X},t_X^*\}=t_X^*$ then $\tilde\gamma[\tilde\beta]=\gamma[\CorrRR{\tilde\beta}]$ and $X(t)\in\overline\Omega_X$. If instead $\min\{t_{0_X},t_X^*\}=t_{0_X}$ then, for $0\leq t\leq t_{0_X}$, $ X(t)\in\overline\Omega_X$. We consider only the case (the other one being easier) $t_{0_X}+k_X\varepsilon_X\leq t\leq t_X^*$. Since $X(t_{0_X})=\tilde X(t_{0_X})\in\partial\Omega_X$ and so $X(t_{0_X})\cdot\xi=\tilde X(t_{0_X})\cdot\xi=0$, we have
\begin{equation}\label{integrals}
\begin{array}{ll}
\displaystyle
\xi\cdot\tilde X(t) =\int_{t_{0_X}}^{t_{0_X}+k_X\varepsilon_X}(f(\tilde X(s))+a_0+D\beta(s))\cdot \xi\ ds \\
\displaystyle
+\int_{t_{0_X}+k_X\varepsilon_X}^t(f(\tilde X(s))+\gamma[\CorrRR{\tilde\beta}](s-k_X\varepsilon_X)+D\beta(s))\cdot \xi\ ds.
\end{array}
\end{equation}
We estimate the first integral in (\ref{integrals}) using (\ref{fscalarexi}). \CorrRR{Indeed, we first assume $t^*_X$ small enough such that $\|\tilde X(s)-X(t_{0_X})\|<r$ as in (\ref{fscalarexi}), for all $s\in[t_{0_X},t^*_X]$, which is possible, independently on the controls and on the points in $K$, because the dynamics are bounded. Here $\zeta$ and $r$ are as in (\ref{fscalarexi}) with respect to $K'$, which is a compact set such that any trajectory starting from a point of $K^X$ does not exit from $K'$ in the time interval $[0,T]$.}
\begin{equation}
\begin{array}{ll}
\label{eq:almostlast}
\displaystyle
\int_{t_{0_X}}^{t_{0_X}+k_X\varepsilon_X}(f(\tilde X(s))+a_0+D\beta(s))\cdot \xi\ ds \geq\\
\displaystyle\zeta k_X\varepsilon_X +\int_{t_{0_X}}^{t_{0_X}+k_X\varepsilon_X}D\beta(s)\cdot \xi \ ds 
\end{array}
\end{equation}

The second integral in (\ref{integrals}) is estimated as
\begin{equation}
\label{eq:last}
\begin{array}{ll}
\displaystyle
\int_{t_{0_X}+k_X\varepsilon_X}^t(f(\tilde X(s))+\gamma[\CorrRR{\tilde\beta}](s-k_X\varepsilon_X)+D\beta(s))\cdot \xi\ ds\\
\displaystyle
= \int_{t_{0_X}}^{t-k_X\varepsilon_X}\bigg(f(\tilde X(s+k_X\varepsilon_X))+\gamma[\CorrRR{\tilde\beta}](s)\\
\displaystyle
+D\beta(s+k_X\varepsilon_X)\bigg)\cdot \xi\ ds \pm \int_{t_{0_X}}^{t-k_X\varepsilon_X}(f(X(s)+D\beta(s)))\cdot \xi \ ds\\
\displaystyle
=\int_{t_{0_X}}^{t-k_X\varepsilon_X}(f(X(s))+\gamma[\CorrRR{\tilde\beta}](s)+D\beta(s))\cdot \xi\ ds\\
\displaystyle
 + \int_{t_{0_X}}^{t-k_X\varepsilon_X} (f(\tilde X(s+k_X\varepsilon_X))-f(X(s))+D\beta(s+k_X\varepsilon_X)-D\beta(s))\cdot\xi \ ds\\
\displaystyle
 \ge(X (t-k_X\varepsilon_X)-X(t_{0_X}))\cdot \xi-Mk_X\varepsilon_X(e^{Lt^*_X} -1)+\\
\displaystyle
\int_{t_{0_X}}^{t-k_X\varepsilon_X} (D\beta(s+k_X\varepsilon_X)-D\beta(s))\cdot\xi \ ds\\
\end{array}
\end{equation}
where in the last inequality (for suitable $M,L>0$ depending on $f$, $A$, $B$ and $D$) we have used the Lipschitz continuity of $f$, standard estimates on trajectories (coming from Gronwall inequality), and the fact that the dynamics is affine in the controls (\ref{nuovaf}). \CorrB{In particular we have used the following estimate for $s\ge t_{0_X}$

\[
\begin{array}{ll}
\displaystyle
|\tilde X(s+k_X\varepsilon_X)-X(s)|\le |\tilde X(t_{0_X}+k_X\varepsilon_X)-X(t_{0_X})|+\\
\displaystyle
\left|\int_{t_{0_X}}^s\left(f(\tilde X(\tau+k_X\varepsilon_X)-f(X(\tau))+\tilde\gamma[\CorrRR{\tilde\beta}](\tau+k_X\varepsilon_X)-\gamma[\CorrRR{\tilde\beta}](\tau)\right.\right.+\\
\displaystyle
D\beta(\tau+k_X\varepsilon_X)-D\beta(\tau)\Big)\Big|d\tau\le\\
\displaystyle
Mk_X\varepsilon_X+L\int_{t_{0_X}}^s|\tilde X(\tau+k_X\varepsilon_X)-X(\tau)|d\tau,
\end{array}
\]

\noindent
using also the equality

\begin{equation}
\label{eq:equality}
\begin{array}{ll}
\displaystyle
\int_{t_{0_X}}^s(D\beta(\tau+k_X\varepsilon_X)-D\beta(\tau))d\tau\\
\displaystyle
=\int_{s-k_X\varepsilon_X}^sD\beta(\tau+k_X\varepsilon_X)d\tau-\int_{t_{0_X}}^{t_{0_X}+k_X\varepsilon_X}D\beta(\tau)d\tau.
\}
\end{array}
\end{equation}

}

Adding (\ref{eq:almostlast}) to (\ref{eq:last}), and using the definition of $\varepsilon_X$, we get
\begin{equation}
\begin{array}{ll}
\displaystyle
\xi\cdot \tilde X(t) \geq \zeta k_X\varepsilon_X -\varepsilon_X-Mk_X\varepsilon_X (e^{Lt_X^*} -1)+ \int_{t-k_X\varepsilon_X}^tD\beta(s)\cdot\xi ds\ge\\
\displaystyle
 (\zeta-\tilde C-M(e^{Lt^*_X}-1))k_X\varepsilon_X-\varepsilon_X
\end{array}
\end{equation}
where $\tilde C$ is an upper bound for $Db\cdot\xi$ and $\zeta-\tilde C>0$ by (\ref{fscalarexi}).
Consequently, if $t^*_X$ is sufficiently small, $\tilde X(t)\cdot \xi \geq \displaystyle\frac{(\zeta-\tilde C) k_X\varepsilon_X}{2}-\varepsilon_X$, and hence, taking $k_X:=2/(\zeta-\tilde C)$, we obtain $\xi\cdot \tilde X(t)\geq 0$. This proves ii), and iv) is proven in a standard way using (\ref{epsilonXestimate}). 

\CorrRR{

\begin{Remark}
\label{rmrk:gluing_decoupling}
Note that $t^*_X$ (as well as $t^*_Y$) only depends, besides the dynamics, on $T>0$ and on the compact set $K$ (via the compact set $K'$ such that, any trajectory starting from $K$, does not exit from $K'$ in the time interval $[0,T]$) and not on controls and non-anticipating strategy in use as well as not on the chosen initial points inside $K$. Hence, repeating, if necessary, the procedure a finite number of time, we can cover the whole interval $[0,\tilde\tau]\subseteq[0,T]$ and obtain the estimates in Assumption 2, where the modulus of continuity ${\cal O}_{T,K}$ only depends on $T$ and $K$. For example, after the first interval $[0,t^*_X]$, we can consider the points $X(t^*_X)$ and $\tilde X(t^*_X)$ and repeat all the construction with those points as initial points (and note that, starting from them, in the time interval $[t^*_X,T]$ we do not exit from $K'$).

We point out once again that the problem is an exit-time problem and hence, when a player firstly exits from its domain, the game stops and what happens after that moment is not influencing anymore. This is the reason why the constructed non-anticipating strategies, even if they map controls on $[0,+\infty[$ to controls on $[0,+\infty[$, are mainly constructed looking to what happens up to the exit time only. Indeed,
after the exit time, controls and strategies may be arbitrarily defined, for example in any a-priori constant manner (which is obviously non-anticipating). Indeed the cost $J$ does not change if we take controls that coincide up to the exit time and possibly differ from the exit time on, since it only depends on controls and strategies used up to exit time.
In particular, a similar situation is in the definitions (\ref{eq:betasegnato}) and (\ref{eq:gammasegnato}), where, after the exit time (of the homologous trajectory), the controls are defined in a constant manner (in that case, in a suitable constant manner, using the outward-pointing control).

In the case when the dynamics $f$ is also strongly decoupled, as in (\ref{eq:hypotheses}) and in the rest of the paper, the definition of $\varepsilon_X$ (\ref{defepsilonX}) is simply replaced by

\[
\varepsilon_X=\sup_{\gamma\in \Gamma}\sup_{\beta\in{\cal B}}\left(\sup_{0\le t\le\min\{\tau_X(x_2, \gamma[\beta]),t_X^*\}}\left(-\xi\cdot  X(t;x_1,\gamma[\beta])\right)^+\right),
\]

\noindent
the definition of $\tilde{\tilde\gamma}$ (\ref{gammatilde}) is replaced by

\[
\begin{array}{l}
\displaystyle
\tilde\gamma[\beta](t)=\left\{
\begin{array}{ll}
\displaystyle
\gamma[\beta](t)&\mbox{if } 0\le t\le\min\{t_{0_X},t_X^*\},\\
\displaystyle
a_0&\mbox{if } \min\{t_{0_X},t_X^*\}\le t\le\min\{t_{0_X},t_X^*\}+k_X\varepsilon_X,\\
\displaystyle
\gamma[\beta](t-k_X\varepsilon_X)&\mbox{if } t\ge\min\{t_{0_X},t_X^*\}+k_X\varepsilon_X,
\end{array}
\right.
\end{array}
\]

\noindent
in place of (\ref{eq:newassumption}) we maintain the corresponding ones ii) and iv) of Assumption 2, and finally (\ref{epsilonXestimate}) turns out as depending only on $\|x_1-x_2\|$.

We finally point out the obvious fact that the estimates iv) and v) of Assumption 2, as well as iv') in (\ref{eq:newassumption}), also hold for all $0\le t\le\tilde\tau$.

\end{Remark}

}

vi) 
\CorrRR{Now we assume that the dependence of the running cost $\ell$ on the controls is separated, that is

\begin{equation}
\label{sumrunningcosts}
\ell(x,y,a,b)=\ell_1(x,y,a)+\ell_2(x,y,b),\ \forall (x,y,a,b)\in\overline\Omega_X\times\overline\Omega_Y\times A\times B,
\end{equation}

\noindent
where $\ell_1,\ell_2$ are continuous, bounded and Lispschitz continuous in $(x,y)$ uniformly with respect to $a$ and $b$ respectively. Note that a similar separated feature is also assumed in Bettiol-Cardaliaguet-Quincampoix \cite{betcarqui}. 
}

We have to estimate
\begin{equation}
\label{costs}
|J_{\tilde\tau}(x_1,y_1,\tilde\gamma[\tilde\beta],\beta)-J_{\tilde\tau}(x_2,y_2,\gamma[\tilde\beta],\tilde\beta)|
\end{equation}
where $\tilde\tau=\min(\tau_X(x_2,\gamma[\tilde\beta],\tilde\beta),\tau_Y(y_1,\beta),T)$ and $J_{\tilde\tau}$ is the integral of the discounted running cost up to time $\tilde\tau$. We will sketch the computation in the case in which $\tilde\tau=T\ge t^*_Y\geq t_{0_X}+k_X\varepsilon_X\geq t_{0_Y}+k_Y\varepsilon_Y\geq t_{0_X}\geq t_{0_Y}$, since the other cases are similar.
The quantity in (\ref{costs}) is majorized by
\begin{equation}
\label{eq:majorization}
\begin{array}{ll}
\displaystyle
 \int_0^{t_{0_Y} }\left|\ell(X(t;x_1),Y(t;y_1),\gamma[\tilde\beta](t),\beta(t))-\ell(X(t;x_2),Y(t;y_2),\gamma[\tilde\beta](t),\beta(t))\right|\,dt\\
\displaystyle
+\int_{t_{0_Y}}^{t_{0_X}}\left|\ell(X(t;x_1),Y(t;y_1),\gamma[\tilde\beta](t),\beta(t))-\ell(X(t;x_2),Y(t;y_2),\gamma[\tilde\beta](t),b_0)\right|\,dt\\
\displaystyle
+\int_{t_{0_X}}^{t_{0_Y}+k_Y\varepsilon_Y}\left|\ell(X(t;x_1),Y(t;y_1),a_0,\beta(t))-\ell(X(t;x_2),Y(t;y_2),\gamma[\tilde\beta](t),b_0)\right|\,dt\\
\displaystyle
+\int_{t_{0_Y}+k_Y\varepsilon_Y}^{t_{0_X}+k_X\varepsilon_X}\Big|\ell(X(t;x_1),Y(t;y_1),a_0,\beta(t))-\\
\displaystyle
\ell(X(t;x_2),Y(t;y_2),\gamma[\tilde\beta](t),\beta(t-k_Y\varepsilon_Y))\Big)\,dt\Big|\\
\displaystyle
+\left|\int_{t_{0_X}+k_X\varepsilon_X}^{t^*_Y}e^{-\lambda t}\Big(\ell(X(t;x_1),Y(t;y_1),\gamma[\tilde\beta](t-k_X\varepsilon_X),\beta(t))-\right.\\
\displaystyle
\left.\phantom{\int_{t_{0_X}+k_X\varepsilon_X}^{t^*_Y}}\ell(X(t;x_2),Y(t;y_2),\gamma[\tilde\beta](t),\beta(t-k_Y\varepsilon_Y))\Big)\,dt\right|
\end{array}
\end{equation}
where we used the expressions of $\tilde\beta$ (\ref{betatilde}) and $\tilde\gamma[\tilde\beta]$ (see (\ref{gammatilde}) and four lines below it)) and dropped the notations of the controls in the trajectories $X(\cdot;x_1)=X(\cdot;x_1,\tilde\gamma[\tilde\beta],\beta)$, $Y(\cdot;y_1)=Y(\cdot;y_1,\beta)$, $X(\cdot;x_2)=X(\cdot;x_2,\gamma[\tilde\beta],\tilde\beta)$, $Y(\cdot;y_2)=Y(\cdot;y_2,\tilde\beta)$. Note that we passed the absolute value under the integral-sign in the first four integrals only.
We are going to use the boundedness and the Lispchitz continuity of $\ell$ (\CorrRR{and of $\ell_1,\ell_2$}), the estimates (\ref{eq:epsilonestimate}), (\ref{epsilonXestimate}), and standard estimates on trajectories.

The first integral in (\ref{eq:majorization}) is majorized by standard procedure \CorrRR{(the difference of the running costs inside the integral is evaluated in the same control values). The other three integrals can be easily estimated since the time interval size is small (it is less than $k_X\varepsilon_X$ or than $k_Y\varepsilon_Y$), while, using (\ref{sumrunningcosts}), the last integral  is majorized by
\begin{equation}\label{runningdecomposed}
\begin{array}{ll}
\displaystyle
\left|\int_{t_{0_X}+k_X\varepsilon_X}^{t^*_Y}e^{-\lambda t}\Big(\ell_1(X(t;x_1),Y(t;y_1),\gamma[\tilde\beta](t-k_X\varepsilon_X))\right.\\
\displaystyle
\left.-\ell_1(X(t;x_2),Y(t;y_2),\gamma[\tilde\beta](t)])\Big)dt\right|\\
\displaystyle
+\left|\int_{t_{0_X}+k_X\varepsilon_X}^{t^*_Y}e^{-\lambda t}\Big(\ell_2(X(t;x_1),Y(t;y_1),\beta(t))\right.\\
\displaystyle
\left.-\ell_2(X(t;x_2),Y(t;y_2),\beta(t-k_Y\varepsilon_Y))\Big)dt\right|.
\end{array}
\end{equation}

\noindent
Now, using the fact that the dynamics are bounded, the fact that the running costs and the function $t\mapsto e^{-\lambda t}$ are bounded and Lipschitz continuous, and the estimates v) of Assumption 2 and iv') of (\ref{eq:newassumption}) (see also Remark \ref{rmrk:gluing_decoupling}),  arguing by a change of variable as in the previous subsection (see for example (\ref{eq:betaestimate}), (\ref{eq:equality}), and also see Bardi-Capuzzo Dolcetta \cite{barcap} page 274), the first integral in (\ref{runningdecomposed}) is estimated by (here $\tilde M$ depends on $L$, the Lipschitz constant of costs and dynamics, on $M$, the bound of costs and dynamics, and on $\lambda$, the discount factor)
\begin{equation}
\begin{array}{ll}
\displaystyle
\displaystyle
\left|\int_{t_{0_X}+k_X\varepsilon_X}^{t^*_Y}\Big(e^{-\lambda(t-k_X\varepsilon_X)}\ell_1(X(t-k_X\varepsilon_X;x_1),Y(t-k_X\varepsilon_X;y_1),\gamma[\tilde\beta](t-k_X\varepsilon_X))\right.\\
\displaystyle
\left.-e^{-\lambda t}\ell_1(X(t;x_2),Y(t;y_2),\gamma[\tilde\beta](t)])\Big)\,dt\right|+\tilde MTk_X\varepsilon_X\le\\
\displaystyle
\Big(\mbox{where we have approximated }  e^{-\lambda t}\ell_1(X(t;x_1),Y(t;y_1),\gamma[\tilde\beta](t-k_x\varepsilon_X)) \mbox{ by}\\
	\displaystyle
e^{-\lambda(t-k_X\varepsilon_X)}\ell_1(X(t-k_X\varepsilon_X;x_1),Y(t-k_X\varepsilon_X;y_1),\gamma[\tilde\beta](t-k_X\varepsilon_X)). \ \mbox{Now, we}\\ \mbox{change the variable in the first addendum inside the integral: $t=t-k_X\varepsilon_X$}\\
\displaystyle
\mbox{and then we pass the absolute value under the integral-sign}\Big)\\
\displaystyle
\le\int_{t_{0_X}}^{t_{0_X}+k_X\varepsilon_X}|\ell_1(X(t;x_1),Y(t;y_1),\gamma[\tilde\beta](t))|dt\\
\displaystyle
+\int_{t_{0_X}+k_X\varepsilon_X}^{t^*_Y-k_X\varepsilon_X}\Big|\ell_1(X(t;x_1),Y(t;y_1),\gamma[\tilde\beta](t))-
\ell_1(X(t;x_2),Y(t;y_2),\gamma[\tilde\beta](t))\Big|\,dt +\\
\displaystyle
\int_{t^*_Y-k_X\varepsilon_X}^{t^*_Y}\Big|\ell_1(X(t;x_2),Y(t;y_2),\gamma[\tilde\beta](t))\Big|\, dt+\tilde MTk_X\varepsilon_X\\\\
\displaystyle
\le Mk_X\varepsilon_X+T{\cal O}_{T,K}(\|x_1-x_2\|+\|y_1-y_2\|)+Mk_X\varepsilon_X+\tilde MTk_X\varepsilon_X\\
\displaystyle
\le {\cal O}_{T,K}(\|x_1-x_2\|+\|y_1-y_2\|),
\end{array}
\end{equation}

\noindent
where in the last inequality we have used (\ref{epsilonXestimate}), and the last ${\cal O}_{T,K}$ is an infinitesimal function, sum of the infinitesimal functions in the line before. We similarly estimate the second integral in (\ref{runningdecomposed}). 

Again, the estimate only depends on $T$ and $K$, being independent on controls, strategies and starting points in $K$. Hence, possibly repeating such procedure a finite number of times, we get the estimate vi) of Assumption 2.

Finally, we point out  that hypothesis (\ref{sumrunningcosts}) is used, here above, in order to separately treat the change of variable in the controls. The extension to the case of non-decoupled cost $\ell(x,y,a,b)$ seems to be not obvious. However, we think that (\ref{sumrunningcosts}) may be probably amended by suitably modifying  (\ref{betatilde}) and (\ref{gammatilde}), but we did not check any details.

}

\begin{acknowledgements}
This work was partially supported by the Italian INDAM-GNAMPA project 2017.
\end{acknowledgements}


\begin{thebibliography}{}

\bibitem{bag} F. Bagagiolo: {\it Minimum time for a hybrid system with thermostatic switchings} in: A. Bemporad, A. Bicchi, G. Buttazzo (Eds.), Hybrid Systems:
Computation and Control, Lecture Notes in Computer Sciences, vol. 4416, Springer-Verlag, Berlin, 2007, pp. 32--45.

\bibitem{bagbar} F. Bagagiolo, M. Bardi: {\it Singular perturbation of a finite horizon problem with state-space constraints}, SIAM J. Control Optim., 36 (1998), 2040--2060.

\bibitem{bagdan} F. Bagagiolo, K. Danieli: {\it Infinite horizon optimal control problems with multiple thermostatic
hybrid dynamics}, Nonlinear Analysis: Hybrid Systems,  6 (2012), 824--838.

\bibitem{bagmag} F. Bagagiolo, R. Maggistro: {\it Hybrid thermostatic approximations
of junctions for some optimal control problems on networks}, SIAM J. Control Optim., 57 (2019), 2415-2442.

\bibitem{barcap} M. Bardi, I. Capuzzo Dolcetta: Optimal Control and Viscosity Solutions of Hamilton-Jacobi-Bellman Equations, Birkh\"auser, Boston 1997.

\bibitem{barkoisor} M. Bardi, S. Koike, P. Soravia: {\it Pursuit-evasion games with state constraints: dynamic programming
and discrete-time approximations}, Discrete Contin. Dyn. Syst., 6 (2000), 361--380.

\bibitem{betcarqui} P. Bettiol, P. Cardaliaguet, M. Quincampoix: {\it Zero-sum state constraint differential game: existence of a value for Bolza problem}, Int. J. Game Theory, 34 (2006), 495--527.

\bibitem{betbrevin} P. Bettiol, A. Bressan, R. Vinter: {\it Estimates for trajectories confined to a cone in $\mathbb{R}^n$}. SIAM J. Control Optim., 49 (2011), 21--41. 

\bibitem{betbreal} P. Bettiol, A. Bressan, R. Vinter: {\it On trajectories satisfying a state constraint: $W^{1,1}$ estimates and counterexamples}, SIAM J. Control Optim., 48 (2010), 4664--4679.

\bibitem{betfac}  P. Bettiol, G. Facchi: {\it Linear estimates for trajectories of state-constrained differential inclusions and normality conditions in optimal control}, J. Math. Anal. Appl. 414 (2014), 914--933.

\bibitem{betfravin}  P. Bettiol, H. Frankowska, R. Vinter: {\it Improved sensitivity relations in state constrained optimal control}, Appl. Math. Optim. 71 (2015), 353--377. 

\bibitem{breal} A. Bressan, G. Facchi: {\it Trajectories of differential inclusions with state constraints}, J. Differential Equations, 250 (2011), 2267--2281.

\bibitem{buccarqui} R. Buckdahn,  P. Cardaliaguet,  M. Quincampoix: {\it Some recent aspects of differential game theory}, Dyn. Games Appl., 1 (2011), 74--114.

\bibitem{carquisain} P. Cardaliaguet, M. Quincampoix, P. Saint-Pierre:  {\it Pursuit differential games with state constraints},
SIAM J. Control Optim. 39 (2000), 1615--1632

\bibitem{dupmce} P. Dupuis, W.M. McEneaney: {\it Risk-sensitive and robust escape criteria}, SIAM J. Control Optim., 35 (1997), 2021--2049.

\bibitem{ellkal} R.J. Elliot, N.J. Kalton: {\it The existence of value in differential games}, Mem. Amer. Math. Soc. 126. AMS, Providence,
USA (1972).

\bibitem{fakoulvin} P. Falugi, P.A. Kountouriotis, R. Vinter: {\it Differential games controllers that confine a system to a safe region in the state space, with applications to surge tank control}, IEEE Trans. Automat. Control, 57 (2012), 2778--2788.

\bibitem{framarmaz}  H. Frankowska, E. M. Marchini, M. Mazzola:  {\it A relaxation result for state constrained inclusions in infinite dimension}, Math. Control Relat. Fields 6 (2016), 113--141.

\bibitem{grogro} D. Gromov, E. Gromova: {\it On a class of hybrid differential games}, Dyn. Games Appl., 7 (2016), 266--288.

\bibitem{koi} S. Koike: {\it On the state constraint problem for differential games},  Indiana Univ. Math. J., 44 (1995), 467--487

\bibitem{son} H. M. Soner: {\it Optimal control problems with state-space constraints I}, SIAM J. Control Optim., 24  (1986), 552--561

\bibitem{Vinter_tank} R. B. Vinter, J. M. C. Clark: {\it A differential games approach to flow control} Proceedings of the 42nd IEEE Conference on Decision and Control Maui, Hawaii USA, December 2003
\end{thebibliography}


\end{document}